\documentclass[11pt]{article}

\usepackage[T1]{fontenc}
\usepackage[french]{babel}
\usepackage{a4wide,epsfig}
\usepackage{amsfonts}
\usepackage{graphicx}
\usepackage{tabularx}
\usepackage{verbatim}
\usepackage{rotate}

\begin{document}

\begin{titlepage}
\begin{center}
\begin{figure}[ht!]
\centerline{
\psfig{figure=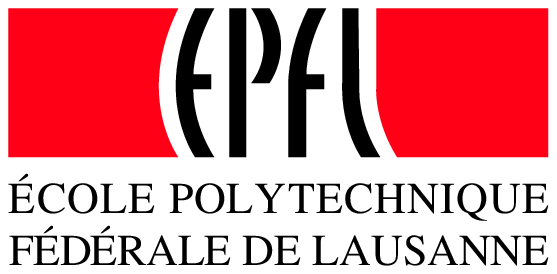,height=3cm}}
\end{figure}
\vspace{0.9cm}
{\LARGE {\underline{Projet STS}}}\\
\vspace{1cm}
{\Huge {\bf Euler et les carrés magiques}}\\
\vspace{0.9cm}  
\psfig{figure=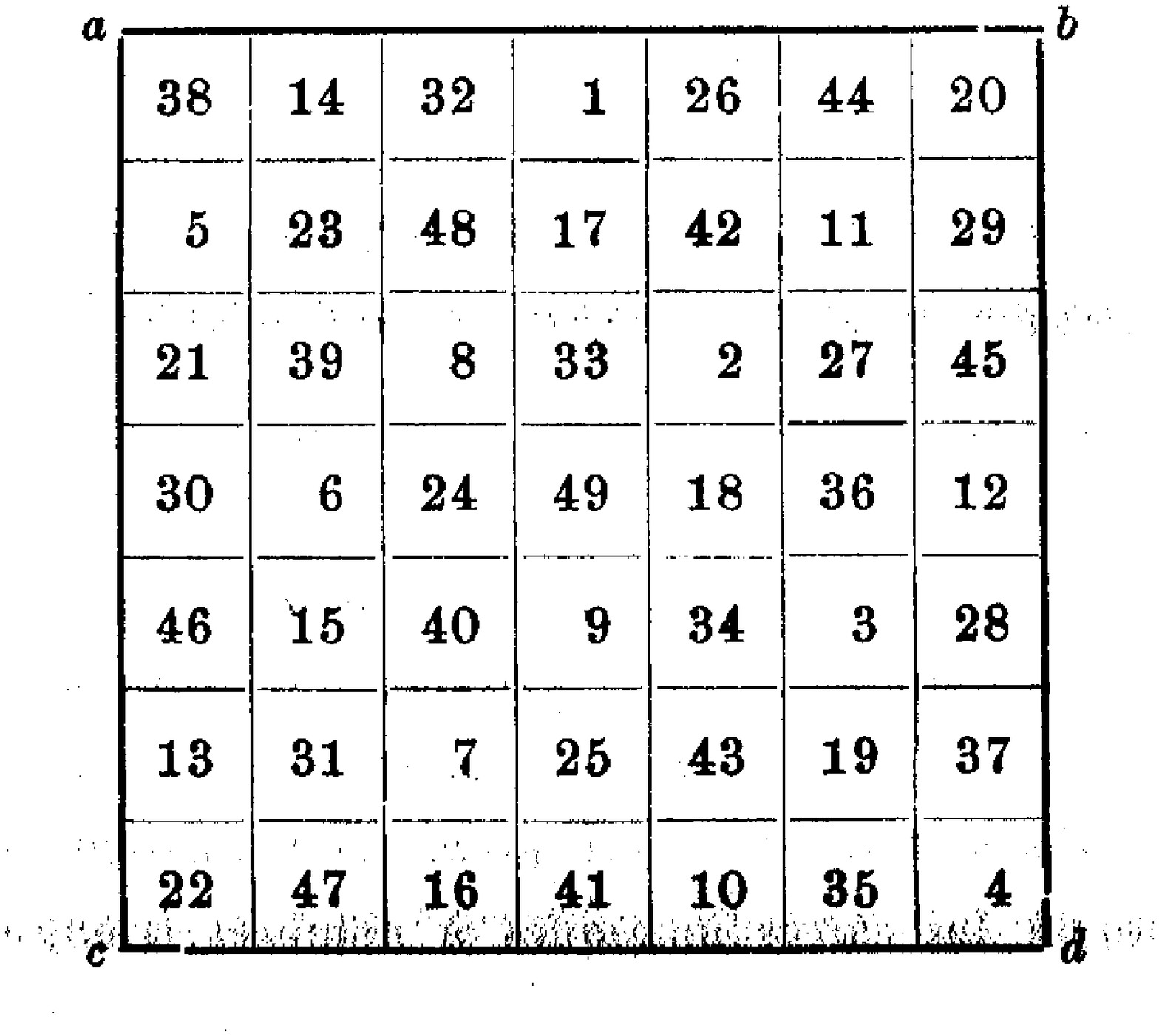,height=5cm}

\vspace{0.9cm}      
{\LARGE  Christophe {\sc Hebeisen}} \\
\vspace{1cm}
{\Large -MA- 4\ ème année}\\
\vspace{1.5cm}
{\LARGE {\underline{Responsable} :}}\\
\vspace{1cm}
{\LARGE Professeur Jacques \textsc{Sesiano}}\\
\vspace{1cm}
{\Large Histoire des mathématiques}\\
\vspace{1cm}
{\Large Printemps 2000-2001}
\end{center}
\end{titlepage}

\newtheorem{thm}{Théorème}
\newtheorem{defin}[thm]{Définition}
\newtheorem{lemme}[thm]{Lemme}
\newtheorem{exple}[thm]{Exemple}
\newtheorem{regle}{Règle}

\tableofcontents
\clearpage

\section{Introduction}
\begin{quotation}
Les carrés magiques ont toujours fasciné la plupart des gens, tant par leur
apparente simplicité que par leur étonnante propriété. Leur origine
est toutefois assez lointaine et incertaine : il n'y a pas de traces de carrés
magiques en Grèce et on trouve seulement un carré de $3\times 3$ en Chine 
vers le début de notre ère.

En Occident, les premiers écrits relatifs aux carrés magiques se trouvent
chez les Arabes%
\footnote{Le premier traité sur ce sujet est attribué au mathématicien
arabe Tâbit ben Korrah.}, au IX$^{\grave{e}me}$ siècle. C'est au
départ une science purement mathémati\-que.

Leur appellation de carré magique provient de leur usage en astrologie et 
comme talismans%
\footnote{On retrouve en Europe deux séries de sept carrés magiques
associés aux planètes et datant du XIV$^{\grave{e}me}$ siècle. \textsc{Cardan} et \textsc{Fermat} s'en sont occupé.}.
On retrouve des carrés magiques sur de nombreuses amulettes au 
XVI$^{\grave{e}me}$ et surtout au XVII$^{\grave{e}me}$ siècle.

\vspace{0.3cm}

{\sc Euler} a consacré deux mémoires et de 
nombreuses pages de ses carnets à l'étude des carrés magiques. 

Dans son mémoire 795, écrit en latin et intitulé \textit{De quadratis 
magicis}, il présente des règles de construction simples pour construire des
carrés magiques à l'aide de lettres latines et grecques auxquelles sont
attribuées des valeurs en progression arithmétique.

Ce mémoire fut présenté à l'Académie des sciences de St-Pétersbourg
le 17 octobre 1776, mais publié seulement en 1849 dans les \textit{Commentationes
arithmeticae}, puis réimprimé en 1862 dans les \textit{Opera postuma}%
\footnote{L'imprimé contient quelques erreurs quant aux valeurs des lettres. Voir traduction.}.

Quant au mémoire 530, intitulé \textit{Recherches sur une nouvelle espèce de
quarrés magiques} et dont le point de départ est le fameux problème des
36 officiers, il fut présenté à cette même académie le 8 mars 1779 et 
fut publié en 1782 dans 
les \textit{Mémoires de la Société des Sciences de Flessingue}, puis
réimprimé en 1849 dans les \textit{Commentationes arithmeticae}. 

\vspace{0.3cm}

La première partie de ce projet est principalement la traduction du texte
\textit{De quadratis magicis} en français.

La deuxième partie est consacrée à l'étude du problème des 36
officiers, qui débouche rapidement sur des considérations sur les carrés
magiques, sous l'aspect des carrés latins à simple marche essentiellement.

\end{quotation}

\pagebreak

\section{\textit{De quadratis magicis} : Sur les carrés magiques}
\label{traduction}
\begin{quotation}

1. Il est coutume d'appeler magique, un carré dont on remplit les cases par
des nombres naturels de telle façon que les sommes des nombres dans chaque
rangée%
\footnote{{\sc Euler} utilise ici le mot \textit{fascia}, qui désigne
une bande ou un ruban ; il sera toujours traduit par <<rangée>>, en
considérant que ce terme désigne autant les lignes et les colonnes que les
deux diagonales du carré en question.}, tant horizontale que
verticale, ainsi que dans les deux diagonales%
\footnote{Sauf mention explicite, 
les diagonales sont ici toujours les diagonales principales.}, soient égales
(entre el\-les)%
\footnote{J'ai essayé de rester aussi fidèle au texte que
possible ; dès lors, je mentionnerai entre parenthèses ( ) les mots du texte 
me semblant superflus, et entre crochets [ ] les mots qui n'apparaissent pas
forcément dans le texte, mais qui à mon avis le complètent ou le rendent 
plus compréhensif.} ; ainsi, si l'on divise les côtés du
carré en $x$ parties égales, le nombre de toutes les cases sera%
\footnote{{\sc Euler} écrit curieusement dans ses textes $xx$ pour désigner $x^2$.} 
$x^2$,
et chaque ligne et colonne, de même que les deux diagonales, auront chacune 
$x$ cases, dans lesquelles il faut donc disposer tous les nombres naturels 1,
2, 3, 4, ... $x^2$ en telle sorte que les sommes pour toutes les ran\-gées
finissent par devenir égales (entre elles). Dans ces circonstances, puisque la
somme de tous ces nombres, de 1 jusqu'à $x^2$, est $$\frac{x^2(1+x^2)}{2}\ ,$$
la somme d'une seule des rangées sera$$=\frac{x(1+x^2)}{2}\ ,$$d'où, si nous
avons $x=3$, la somme pour une seule rangée sera $=15$.

\vspace{0.3cm}
2. Donc, à partir de ceci, quel que soit le nombre de cases dans lequel le
carré tout entier est divisé, on pourra facilement déterminer la somme des
nombres disposés dans une rangée, et il nous sera utile que cette table ait
été établie pour chaque rangée pour tout carré de ce genre :
$$
\begin{array}{c|c|c}
\vspace{0.1cm}
x & x^2 & \frac{x(1+x^2)}{2} \\
\hline 
1 & 1   & 1   \\
2 & 4   & 5   \\
3 & 9   & 15  \\
4 & 16  & 34  \\ 
5 & 25  & 65  \\
6 & 36  & 111 \\
7 & 49  & 175 \\
8 & 64  & 260 \\
9 & 81  & \,\,369\footnotemark\\
  & etc.&
\end{array}
$$%
\footnotetext{360 dans les \textit{Opera postuma}.} 
où $x$ désigne le nombre de parties par lesquelles le côté du carré est
divisé, $x^2$ le nombre de cases contenues dans le carré et
$\frac{1}{2}x(1+x^2)$ indique la somme de tous les nombres disposés dans une
rangée.

\vspace{0.3cm}
3. Afin de rechercher une règle précise pour construire de tels carrés 
magiques de n'importe quel ordre, il est de la plus grande
importance de remarquer que chaque nombre de 1, 2, 3 etc. jusqu'à $x^2$ peut 
être représenté par la formule $$mx+n\,.$$
En effet, si nous remplaçons $m$ par les valeurs successives
0, 1, 2, 3, 4 jusqu'à $x-1$, puis $n$ par les
valeurs 1, 2, 3, 4, ... $x$, il est par là clair que l'on obtiendra tous les 
nombres de 1 à $x^2$, puisque chaque valeur de $n$ est successivement
combinée avec toutes les valeurs de $m$%
\footnote{Cela n'est pas si évident à première vue ; une représentation 
sous forme de tableau permet de s'en faire une meilleure idée : 
$$
\begin{array}{c|ccccc}
\begin{array}{cc} \ddots & n\\ m & \ddots \end{array} & 1 & 2 & 3 & \cdots & x \\
\hline
0 &   1   &    2   &  3  & \cdots & x  \\
1 & x+1  &   x+2  & x+3 & \cdots & 2x \\
\vdots & \vdots & \vdots & \vdots & \ddots & \vdots \\
x-1 & (x-1)\,x+1  & (x-1)\,x+2 & (x-1)\,x+3 & \cdots & x^2 
\end{array}
$$ 
comme les valeurs maximales que peuvent prendre $m$ et $n$ sont $x-1$ et $x$
respectivement, l'expression ne dépassera donc jamais $(x-1)x+x=x^2$, ce qui 
est bien demandé.}. 
Puisque tous les nombres inscrits dans le carré peuvent être traduits de 
cette manière par la formule $mx+n$, et donc représentés en deux parties, 
nous désignerons par la suite les premières parties, $mx$, simplement par 
les lettres latines $a$, $b$, $c$, $d$, etc. et les deuxièmes, $n$, par 
les lettres grecques \textit{$\alpha$, $\beta$, $\gamma$, $\delta$}, où il est 
clair que pour n'importe quel nombre $x$ la quantité des lettres, tant latines
que grecques, doit être $x$, puisque les valeurs des lettres latines seront 
$0x$, $1x$, $2x$, $3x$ jusqu'à $(x-1)x$, tandis que celles des lettres 
grecques sont 1, 2, 3, 4, ...~$x$. Et il ne faut pas ici s'imaginer devoir 
choisir un ordre déterminé dans ces lettres, tant latines que grecques, 
puisque n'importe quelle lettre latine peut représenter à volonté soit 
$0x$, soit $1x$, soit $2x$ etc., pourvu qu'on attribue des valeurs distinctes
à chaque lettre ; et il faut faire de même pour les lettres grecques.

\vspace{0.3cm}
4. Dans ce qui suit on pourra donc représenter n'importe quel nombre à
inscrire dans le carré par une somme d'une lettre latine et d'une grecque, 
par exemple par
$b+\delta$ ou [encore] $a+\beta$ etc., en telle sorte que chacun des
nombres soit représenté en deux parties ; alors en effet, si l'on réunit
chaque lettre latine avec chaque lettre grecque, clairement il doit en résulter
tous les nombres de 1 jusqu'à $x^2$ ; et de même il est évident que les
diverses combinaisons de ces lettres forment toujours des nombres différents,
et qu'aucun nombre ne peut être exprimé de deux manières [différentes].

\vspace{0.3cm}
5. Donc, puisque tous les nombres sont représentés par la réunion d'une 
lettre latine et grecque, établissons la règle fondamentale suivante pour la
construction des carrés magiques : à savoir [que] d'abord j'inscris les
lettres latines dans chaque case du carré en telle sorte que leur somme soit
la même dans chaque rangée, où, comme le nombre de ces lettres est $=x$,
[et que] d'autre part le nombre de toutes les cases est $=x^2$, il est
évident que l'on doive répéter n'importe quelle lettre $x$ fois. D'autre 
part, on comprend que, de façon identique, on inscrit les lettres grecques
dans les cases du même carré en sorte que leurs sommes dans chaque rangée
deviennent égales. Ainsi, encore une fois, les sommes de tous les nombres
formés d'une lettre latine et grecque seront, pour chaque rangée, égales 
(entre elles). Il ne reste plus que, dans cette disposition, à chaque lettre 
latine différente soit associée une lettre grecque différente,
puisque par ce procédé aucun nombre de 1 à $x^2$ ne sera omis, ni ne
pourra apparaître deux fois.

\vspace{0.3cm}
6. Une fois ces règles établies généralement, 
étudions les différentes formes des carrés suivant le nombre de cases : il
apparaît immédiatement que celui-ci commence à neuf, puisque l'on ne 
peut disposer des nombres de la façon décrite ci-dessus dans le carré 
divisé en quatre cases. De plus, il sera utile de remarquer que, 
indépendamment de l'ordre,
puisque pour n'importe quelle catégorie [de carré] le nombre de lettres tant
latines et grecques est $=x$, et que chaque rangée contient ce même nombre
de cases, les conditions requises seront satisfaites si nous inscrivons dans
une rangée toutes les lettres distinctes, tant latines que grecques. Mais s'il
advient que la même lettre apparaisse deux ou trois fois dans
quelque rangée, il est toujours nécessaire que la somme de toutes
les lettres dans une même rangée soit égale à la somme de toutes les
lettres, soit latines $a+b+c+d\,+$ etc., soit grecques
$\alpha+\beta+\gamma+\delta\,+$ etc.

\end{quotation}

\subsection{La catégorie des carrés divisés en 9 cases}
\label{9}
\begin{quotation}

7. Vu que pour cette espèce [de carré] $x=3$, nous aurons le même nombre 
de lettres latines $a$, $b$, $c$, et grecques $\alpha$, $\beta$,
$\gamma$, les valeurs des lettres latines seront ici 0, 3%
\footnote{2 dans les \textit{Opera postuma}}, 6, et celles des
lettres grecques 1, 2, 3. Commençons maintenant par les lettres latines 
$a$, $b$, $c$ : il sera facile de les inscrire dans notre carré divisé en 9
cases en telle sorte que dans chaque rangée, tant horizontale que verticale,
chacune de ces trois lettres apparaisse [une seule fois], ce que l'on peut voir
dans cet exemple
$$
\begin{array}{ccc}
\vspace{0.1cm}
a & b & c \\
\vspace{0.1cm}
b & c & a \\
c & a & b
\end{array}
$$
où l'on retrouve ces mêmes trois lettres $a$, $b$, $c$ dans l'une des deux
diagonales, alors que dans l'autre la lettre $c$ se répète trois fois ;
on remarque également sans difficulté qu'il ne peut clairement pas
arriver que toutes les lettres des deux diagonales soient simultanément
différentes ;  ce fait ne dérange nullement pour autant que la somme des
lettres de cette diagonale, à savoir $3c$, égale celle des rangées
restantes $a+b+c$ ; c'est-à-dire, pour autant que l'on ait $2c=a+b$. De là
il est clair que $c$ doit prendre la valeur 3, et que l'on doive assigner aux
lettres $a$ et $b$ les valeurs 0 et 6 ; en effet, ainsi nous avons que
$2c=a+b$. Or on peut%
\footnote{{\sc Euler} utilise souvent d'autres formes temporelles (ici p.ex. un
futur), que j'ai souvent traduites par un autre temps, de façon à garder 
une certaine continuité dans le texte. De même j'ai souvent omis certains 
mots dont il abuse (notamment \textit{vero, autem, etc.}) pour accentuer le poids 
de ses affirmations.}
poser à choix soit $a=0$, soit $b=0$ ; observant ceci il résulte que la 
somme de chaque rangée est $a+b+c=9$.

\vspace{0.3cm}
8. Il sera loisible de répartir les lettres grecques dans un tel carré de 
façon semblable ; représentons-les dans la même figure dans l'ordre 
inverse :
$$\begin{tabular}{ccc}
\vspace{0.2cm}
$\gamma$ & $\beta$  & $\alpha$ \\
\vspace{0.2cm}
$\alpha$ & $\gamma$ & $\beta$ \\
$\beta$  & $\alpha$ & $\gamma$
\end{tabular}$$
il est nécessaire d'y avoir $2\gamma=\alpha+\beta$ et par là-même
$\gamma=2$. Ainsi en effet, si nous combinons chaque case de la première 
figure avec chacune de celle-ci dans l'ordre naturel, il est évident que 
n'importe quelle lettre latine sera associée avec chacune des lettres 
grecques, en sorte que de cette conjonction résulteront tous les nombres
de 1 jusqu'à 9 ; cette combinaison produit la figure suivante :
$$\begin{tabular}{ccc}
\vspace{0.4cm}
$a\gamma$ & $b\beta$  & $c\alpha$ \\
\vspace{0.4cm}
$b\alpha$ & $c\gamma$ & $a\beta$ \\
$c\beta$  & $a\alpha$ & $b\gamma$
\end{tabular}$$ 
où l'on notera que les paires de lettres jointes ne représentent pas un
produit, mais une somme.

\vspace{0.3cm}
9. Donc dans cette figure on doit prendre les valeurs $c=3$ et
$\gamma=2$, en sorte que l'on doit attribuer aux lettres $a$ et $b$ les
valeurs 0 et 6 [d'une part], et aux lettres $\alpha$ et $\beta$ les valeurs
1 et 3 [d'autre part] ; si nous posons $a=0$ et $b=6$, puis $\alpha=1$ et
$\beta=3$, il en résulte le carré magique suivant : 
$$\begin{tabular}{lccc}
\vspace{0.1cm}
  & 2 & 9 & 4 \\
\vspace{0.1cm}
I)& 7 & 5 & 3 \\
  & 6 & 1 & 8
\end{tabular}$$
où la somme de n'importe quelle rangée donne 15. Si nous voulons échanger
les valeurs des lettres $a$ et $b$, de même que celles de $\alpha$ et
$\beta$, on remarque facilement que [de cette opération] on changera
seulement la position du carré%
\footnote{En effet, les valeurs de $c$ et $\gamma$ étant fixées ($c=3$,
$\gamma=2$), les trois autres possibilités sont :
$$\begin{tabular}{lcccclcccclccc}
\vspace{0.1cm}
    & 8 & 3 & 4 &&      & 2 & 7 & 6 &&     & 8 & 1 & 6 \\
\vspace{0.1cm}
II) & 1 & 5 & 9 && III) & 9 & 5 & 1 && IV) & 3 & 5 & 7 \\
\vspace{0.2cm} 
    & 6 & 7 & 2 &&      & 4 & 3 & 8 &&     & 4 & 9 & 2 \\
&\multicolumn{3}{c}{a=6, b=0} &&
&\multicolumn{3}{c}{a=0, b=6} &&
&\multicolumn{3}{c}{a=6, b=0} \\
&\multicolumn{3}{c}{$\alpha$=1, $\beta$=3} &&
&\multicolumn{3}{c}{$\alpha$=3, $\beta$=1} &&
&\multicolumn{3}{c}{$\alpha$=3, $\beta$=1} \\
\end{tabular}$$}.

\vspace{0.3cm}
10. Certes cette disposition, tant des lettres latines que grecques, est
suffisamment claire en elle-même, mais l'importance particulière réside en
ceci qu'il a été posé que par la combinaison faite chaque lettre latine 
est associée à chaque lettre grecque%
\footnote{Donc la lettre $a$ est associée avec chacune des trois lettres
grecques $\alpha$, $\beta$, et $\gamma$, et ainsi de suite pour chaque lettre.},
ce qui dans notre disposition apparaît avoir été obtenu par hasard. Or, 
afin de ne rien laisser à l'arbitraire dans notre discussion, remarquons avant
tout que l'ordre des lettres grecques $\alpha$, $\beta$, $\gamma$ ne dépend en
aucun cas de celui des lettres latines $a$, $b$, $c$, en sorte que pour
n'importe quelle rangée définie à l'aide de lettres latines, on aurait
pu choisir de les combiner avec les lettres grecques portant le même nom, à
savoir $\alpha$ avec $a$, $\beta$ avec $b$ et $\gamma$ avec $c$ ;
ainsi, si l'on place dans la première rangée horizontale $a\alpha$,
$b\beta$, $c\gamma$, [et] puisque la même lettre grecque ne doit pas
apparaître deux fois dans quelque rangée horizontale ou verticale, il
est sans conteste évident que la deuxième rangée horizontale sera
$b\gamma$, $c\alpha$, $a\beta$ et la troisième $c\beta$, 
$a\gamma$, $b\alpha$ ; d'où le carré :
$$\begin{tabular}{ccc}
\vspace{0.4cm}
$a\alpha$ & $b\beta$ & $c\gamma$ \\
\vspace{0.4cm}
$b\gamma$ & $c\alpha$ & $a\beta$ \\
$c\beta$ & $a\gamma$ & $b\alpha$
\end{tabular}$$ 
où, du fait que dans la diagonale de gauche la même lettre grecque $\alpha$
apparaît trois fois, il est nécessaire d'avoir
$3\alpha=\alpha+\beta+\gamma$, et à cause de cela $2\alpha=\beta+\gamma$ ;
c'est pourquoi, à partir de là, la valeur même de $\alpha$ est
[entièrement] déterminée, à savoir $\alpha=2$, de même que nous voyons
que $c$ doit prendre la valeur $c=3$ ; mais de là ne naissent pas de nouveaux
carrés magiques%
\footnote{A nouveau, les valeurs de $c$ et $\alpha$ étant fixées ($c=3$,
$\alpha=2$), on retrouve les quatre possibilités d'avant :
$$\begin{tabular}{lcccclcccclcccclccc}
\vspace{0.1cm}
   & 2 & 9 & 4 &&     & 8 & 3 & 4 &&      & 2 & 7 & 6 &&     & 8 & 1 & 6 \\
\vspace{0.1cm}
I) & 7 & 5 & 3 && II) & 1 & 5 & 9 && III) & 9 & 5 & 1 && IV) & 3 & 5 & 7 \\
\vspace{0.2cm} 
   & 6 & 1 & 8 &&     & 6 & 7 & 2 &&      & 4 & 3 & 8 &&     & 4 & 9 & 2 \\    
&\multicolumn{3}{c}{a=0, b=6} &&
&\multicolumn{3}{c}{a=6, b=0} &&
&\multicolumn{3}{c}{a=0, b=6} &&
&\multicolumn{3}{c}{a=6, b=0} \\
&\multicolumn{3}{c}{$\beta$=3, $\gamma$=1} &&
&\multicolumn{3}{c}{$\beta$=3, $\gamma$=1} &&
&\multicolumn{3}{c}{$\beta$=1, $\gamma$=3} &&
&\multicolumn{3}{c}{$\beta$=1, $\gamma$=3} \\
\end{tabular}$$}.

\vspace{0.3cm}
11. Bien que l'arrangement des lettres grecques ne présente aucune
difficulté dans ce premier genre [de carré], pour ceux qui ont un nombre de
cases plus élevé, il importe d'apporter une règle fixe permettant 
d'inscrire les lettres grecques une fois que les lettres latines auront été 
placées de façon convenable ; à cette fin choisissons quelque rangée 
du milieu, soit horizontale, soit verticale, soit même diagonale, telle que de
part et d'autre de cette rangée nous retrouvions dans toutes les cases
équidistantes deux lettres latines différentes. C'est le cas pour la colonne%
\footnote{{\sc Euler} utilise ici le terme \textit{columna}, qui désigne autant
une ligne qu'une colonne, contrairement à notre habitude. Pour préciser, il
rajoute donc colonne <<horizontale>> ou <<verticale>>, que j'ai la plupart du temps
traduit par ligne ou colonne.}
du milieu, autour de laquelle nous retrouvons dans la première ligne les 
lettres $a$ et $c$, dans la deuxième $b$ et $a$ et dans la troisième $c$ 
et $b$, où partout deux lettres différentes se font face%
\footnote{Voir le carré au \S.10 ci-dessus.}.

\vspace{0.3cm}
12. Une fois trouvée une telle rangée médiane, associons à chaque lettre
latine du carré les lettres grecques de même nom, puis échangeons les
lettres grecques dans les cases se correspondant de part et d'autre [de cette
rangée]%
\footnote{{\sc Euler} procède donc de la façon suivante dans son exemple :
il part du carré rempli par les lettres latines, puis repère une rangée
médiane (i) ; il complète avec les lettres grecques (ii), et finalement il 
échange les lettres grecques symétriquement par rapport à cette randée
médiane (iii) :
$$\begin{tabular}{lc|c|cclcccclccc}
\cline{3-3}
\vspace{0.3cm}
    & $a$ & $b$ & $c$ &               &      & $a\alpha$ & $b\beta$  & $c\gamma$ &               &       & $a\gamma$ & $b\beta$  & $c\alpha$ \\
\vspace{0.3cm}    
(i) & $b$ & $c$ & $a$ & $\Rightarrow$ & (ii) & $b\beta$  & $c\gamma$ & $a\alpha$ & $\Rightarrow$ & (iii) & $b\alpha$ & $c\gamma$ & $a\beta$  \\
    & $c$ & $a$ & $b$ &               &      & $c\gamma$ & $a\alpha$ & $b\beta$  &               &       & $c\beta$  & $a\alpha$ & $b\gamma$ \\
\cline{3-3}   
\end{tabular}$$} ; 
de cette façon de faire il résultera cette figure :
$$\begin{tabular}{ccc}
\vspace{0.4cm}
$a\gamma$ & $b\beta$  & $c\alpha$ \\
\vspace{0.4cm}
$b\alpha$ & $c\gamma$ & $a\beta$  \\
$c\beta$  & $a\alpha$ & $b\gamma$ 
\end{tabular}$$
où nous sommes sûrs qu'à chaque lettre latine est associée chaque
lettre grecque.
Pour le reste, afin que la condition des diagonales soit satisfaite, il faut
avoir, comme nous l'avions déjà mentionné, $$2c=a+b\ \ \ 
\mbox{ et } \ \ \ 2\gamma=\alpha+\beta\,.$$
Mais cette figure ne diffère pas de celle que nous avions trouvée au \S.8.
Finalement il faut remarquer que de quelque manière que l'on permute les
rangées horizontales et verticales, la somme de ces rangées ne change pas.
Cependant, dans les diagonales peut survenir une grande différence ; ainsi, si
nous enlevons la première colonne et que nous l'apposons à droite%
\footnote{{\sc Euler} met : à gauche ($!$), ce qui est manifestement une erreur,
si l'on considère l'orientation habituelle de gauche à droite (voir
également la note au \S.17).}
[du reste] apparaît cette figure :
$$\begin{tabular}{ccc}
\vspace{0.4cm}
$b\beta$  & $c\alpha$ & $a\gamma$ \\
\vspace{0.4cm}
$c\gamma$ & $a\beta$  & $b\alpha$ \\
$a\alpha$ & $b\gamma$ & $c\beta$  
\end{tabular}$$
où à cause des diagonales il faut avoir $$2a=b+c\ \ \ \mbox{ et } \ \ \ 
2\beta=\alpha+\gamma\ ;$$  
il faut remarquer ceci pour toutes les transpositions [possibles]. Cette 
observation sera d'une importance capitale pour les catégories suivantes.
\end{quotation}

\subsection{La catégorie des carrés divisés en 16 cases}
\label{16}
\begin{quotation}
13. Puisque pour cette catégorie nous avons $x=4$, nous aurons quatre lettres
latines $a$, $b$, $c$, $d$, dont les valeurs sont 0, 4, 8, 12, et ce même
nombre de lettres grecques $\alpha$, $\beta$, $\gamma$, $\delta$, dont les
valeurs seront 1, 2, 3, 4. En premier lieu, inscrivons dans un tel carré les
quatre lettres latines en telle sorte que dans toutes les rangées, tant 
horizontales que verticales, chacune de ces quatre lettres apparaisse, et en 
sorte que cela se produise également dans les deux diagonales, si cela se 
peut. 

\vspace{0.3cm}
14. Or puisque parmi ces lettres $a$, $b$, $c$, $d$\, aucun ordre n'est prescrit,
inscrivons-les  dans la première ligne, dans l'ordre, et aussi dans la
diagonale de gauche, où dans la deuxième case de celle-ci on pourra
écrire soit la lettre $c$, soit $d$ ; écrivons donc $c$ : désormais toutes
les rangées qui restent sont [entièrement] déterminées, pourvu que l'on
fasse attention de ne pas mettre deux fois la même lettre dans la même
rangée, aussi bien horizontale que verticale%
\footnote{En effet, une fois la première ligne et la deuxième case de la
diagonale remplies (i), nous pouvons compléter la diagonale : comme il nous
reste à disposition les lettres $b$ et $d$, et pour que la
dernière colonne ne contienne pas deux fois la lettre $d$, il nous faut
placer celle-ci dans la troisième case de la diagonale, et $b$
dans la quatrième (ii) ; nous pouvons procéder de même avec la deuxième
colonne (iii), par exemple, et il ne reste plus qu'à remplir le carré...
$$\begin{tabular}{lccccclccccclcccc}
\vspace{0.1cm}
    & $a$     & $b$     & $c$     & $d$     &               &      & $a$     & $b$     & $c$     & $d$     &               &       & $a$     & $b$ & $c$     & $d$     \\
\vspace{0.1cm}
(i) & $\cdot$ & $c$     & $\cdot$ & $\cdot$ & $\Rightarrow$ & (ii) & $\cdot$ & $c$     & $\cdot$ & $\cdot$ & $\Rightarrow$ & (iii) & $\cdot$ & $c$ & $\cdot$ & $\cdot$ \\
\vspace{0.1cm}   
    & $\cdot$ & $\cdot$ & $\cdot$ & $\cdot$ &               &      & $\cdot$ & $\cdot$ & $d$     & $\cdot$ &               &       & $\cdot$ & $a$ & $d$     & $\cdot$ \\
    & $\cdot$ & $\cdot$ & $\cdot$ & $\cdot$ &               &      & $\cdot$ & $\cdot$ & $\cdot$ & $b$     &               &       & $\cdot$ & $d$ & $\cdot$ & $b$  
\end{tabular}$$} ; 
en procédant ainsi nous trouverons la figure suivante :
$$
\begin{array}{cccc}
\vspace{0.1cm}
a & b & c & d \\
\vspace{0.1cm}
d & c & b & a \\
\vspace{0.1cm}
b & a & d & c \\
c & d & a & b
\end{array}
$$
où de surcroît même l'autre diagonale contient toutes les quatre 
lettres, en sorte qu'alors aucune condition n'est prescrite sur les valeurs des 
lettres $a$, $b$, $c$, $d$. Ici nous aurions pu également inscrire la lettre 
$d$ dans la deuxième case de la diagonale, mais la figure résultant [de 
cette opération] ne différerait en aucune autre façon de cette figure, 
si ce n'est par l'emplacement [des lignes]%
\footnote{Effectivement, si nous construisons le carré avec la lettre $d$ dans
la deuxième case de la diagonale, nous obtenons le carré suivant :
$$
\begin{array}{cccc}
\vspace{0.1cm}
a & b & c & d \\
\vspace{0.1cm}
c & d & a & b \\
\vspace{0.1cm}
d & c & b & a \\
b & a & d & c 
\end{array}
$$}, 
en sorte que cette figure doit être considérée [comme celle] qui englobe 
tous les cas possibles.

\vspace{0.3cm}
15. Maintenant pour l'insertion des lettres grecques : puisqu'il n'y a aucune
rangée médiane, ni parmi les horizontales, ni les verticales%
\footnote{Revoir la démarche utilisée pour les carrés d'ordre 3 au \S.11.},
prenons donc la diagonale [contenant] $a$, $c$, $d$, $b$ comme rangée du 
milieu ; il nous apparaîtra bientôt que dans les cases à même 
distance de part et d'autre et se correspondant%
\footnote{Donc symétriques par rapport à cette diagonale.}
se retrouvent deux lettres différentes (entre elles) ;
nous pouvons donc utiliser la règle donnée au \S.11 en toute sécurité. 
En premier lieu donc, associons aux lettres placées dans cette diagonale les
lettres grecques de même nom, ensuite échangeons les lettres grecques dans 
les cases correspondantes ; de cette manière sera formée la figure 
suivante :
$$
\begin{tabular}{cccc}
\vspace{0.4cm}                                
$a\alpha$ & $b\delta$ & $c\beta$  & $d\gamma$ \\
\vspace{0.4cm}                                
$d\beta$  & $c\gamma$ & $b\alpha$ & $a\delta$ \\
\vspace{0.4cm}                                
$b\gamma$ & $a\beta$  & $d\delta$ & $c\alpha$ \\                                
$c\delta$ & $d\alpha$ & $a\gamma$ & $b\beta$
\end{tabular}
$$

\vspace{0.3cm}
16. Dans cette figure donc toutes les quatre lettres, tant latines que grecques,
apparaissent dans chaque rangée, tant horizontale que verticale et diagonale ;
d'où l'on peut attribuer à volonté chaque fois quatre valeurs numériques
à ces lettres, et sans aucune restriction. Or comme de ces quatre lettres on 
peut en tirer 24 arrangements%
\footnote{En effet, le nombre de possibilités de placer quatre lettres
distinctes sans répétition est de $4 !=24$.\label{note}},
on pourra former au total 576 figures différentes%
\footnote{C'est-à-dire $24\cdot 24$, une fois pour chaque alphabet.}, 
où il est certain que plusieurs [d'entre elles] ne différeront l'une de
l'autre seulement par raison de la disposition. 

\vspace{0.3cm}
17. Mais de ce qui précède il ne faut nullement conclure que cette figure
contienne complètement tous les carrés magiques de cette espèce. En effet
outre ceux-ci on peut en trouver où chacune des quatre lettres, tant latines
que grecques, ne se retrouve pas, mais qui ne remplissent pas moins les
conditions prescrites ; or de telles formes peuvent apparaître par 
déplacements de lignes ou de colonnes ; par exemple si dans la figure 
ci-dessus on déplace la première colonne à la fin%
\footnote{Ici l'orientation est bien de gauche à droite.},
il apparaît cette figure :
$$
\begin{tabular}{cccc}
\vspace{0.4cm}                                
$b\delta$ & $c\beta$  & $d\gamma$ & $a\alpha$ \\
\vspace{0.4cm}                                
$c\gamma$ & $b\alpha$ & $a\delta$ & $d\beta$  \\
\vspace{0.4cm}                                
$a\beta$  & $d\delta$ & $c\alpha$ & $b\gamma$ \\                                
$d\alpha$ & $a\gamma$ & $b\beta$  & $c\delta$ 
\end{tabular}
$$
où l'on retrouve certes encore chaque lettre, tant latine que grecque, dans
toutes les lignes et colonnes, mais où dans la diagonale descendant de gauche
à droite on rencontre seulement deux lettres latines, à savoir $b$ et $c$,
et également seulement deux lettres grecques, $\alpha$ et $\delta$. Au
contraire dans l'autre diagonale [on trouve] seulement les deux
lettres latines $a$ et $d$, alors que [l'on retrouve], comme auparavant,
seulement les lettres grecques $\alpha$ et $\delta$.

\vspace{0.3cm}
18. Pour que cette figure satisfasse aux conditions prescrites, on ne peut en
outre attribuer [arbitrairement] des valeurs numériques particulières aux
lettres, mais la condition doit être rajoutée d'avoir, pour les lettres
latines,
$$
b+c=a+d\,,
$$
et pour les lettres grecques
$$
\alpha+\delta=\beta+\gamma\,;
$$
aussi si nous choisissons $a=0$, il faut poser $d=12$, afin d'avoir $b=4$ et
$c=8$, ou vice versa $c=4$ et $b=8$. De la même manière si pour les lettres
grecques nous choisissons $\alpha=1$, il faut avoir $\delta=4$, et alors
$\beta=2$ et $\gamma=3$. De là est formé le carré magique entièrement
déterminé
$$
\begin{tabular}{cccc}
\vspace{0.3cm}                                
$8$  & $10$ & $15$ & $1$  \\
\vspace{0.3cm}                                
$11$ & $5$  & $4$  & $14$ \\
\vspace{0.3cm}                                
$2$  & $16$ & $9$  & $7$  \\                                
$13$ & $3$  & $6$  & $12$ 
\end{tabular}
$$
où il est clair que la somme de chaque rangée est 34. Mais de telles formes
déterminées pourront être formées de nombreuses autres façons par 
des transpositions de lignes ou colonnes.

\vspace{0.3cm}
19. Et il n'est également pas absolument requis que dans chaque ligne ou
colonne apparaissent toutes les lettres, tant latines que grecques ; il se 
peut même que dans ces lignes ou colonnes apparaissent seulement deux lettres
latines ou grecques pourvu que leur somme soit égale à la moitié de
toutes les quatre. Pour de telles figures on a besoin d'opérations
particulières, pour lesquelles des règles fixes peuvent être difficilement
établies ; il faut seulement que les lettres latines et grecques soient 
disposées en telle sorte que non seulement elles fassent la somme due dans 
chaque rangée, mais encore que toutes les lettres grecques soient associées 
avec chaque lettre latine%
\footnote{C'est-à-dire les conditions habituelles des carrés magiques, à
savoir :
\begin{itemize}
    \item des sommes égales dans toutes les lignes, toutes les colonnes et
    les deux diagonales ;
    \item des nombres tous différents et tous utilisés.
\end{itemize}}.

\vspace{0.3cm}
20. Afin de donner un exemple d'une telle opération, posons d'abord que
$$a+d=b+c$$ et disposons les lettres latines comme suit 
$$
\begin{tabular}{cccc}
\vspace{0.1cm}                                
$a$ & $a$ & $d$ & $d$ \\
\vspace{0.1cm}                                
$d$ & $d$ & $a$ & $a$ \\
\vspace{0.1cm}                                
$b$ & $b$ & $c$ & $c$ \\                                
$c$ & $c$ & $b$ & $b$ 
\end{tabular}
$$
où pour chaque rangée la somme des nombres est certainement la même ; 
quant aux lettres grecques, associons aux lettres latines de la diagonale de 
gauche les lettres grecques de même nom, [et] puisque nous constatons que de
part et d'autre de cette diagonale sont disposées chaque fois deux lettres
différentes, décidons de leur associer les lettres grecques permutées, ce
qui donne la figure suivante :
$$
\begin{tabular}{cccc}
\vspace{0.4cm}                                
$a\alpha$ & $a\delta$ & $d\beta$  & $d\gamma$ \\
\vspace{0.4cm}                                
$d\alpha$ & $d\delta$ & $a\beta$  & $a\gamma$ \\
\vspace{0.4cm}                                
$b\delta$ & $b\alpha$ & $c\gamma$ & $c\beta$  \\                         
$c\delta$ & $c\alpha$ & $b\gamma$ & $b\beta$ 
\end{tabular}
$$
où donc il est nécessaire d'avoir pour les lettres grecques
$$\alpha+\delta=\beta+\gamma\,;$$
ainsi si nous prenons $a=0$, $b=4$, $c=8$, $d=12$ et $\alpha=1$, $\beta=2$,
$\gamma=3$ et $\delta=4$, apparaît ce carré magique :
$$
\begin{tabular}{cccc}
\vspace{0.3cm}                                
$1$  & $4$  & $14$ & $15$ \\
\vspace{0.3cm}                                
$13$ & $16$ & $2$  & $3$  \\
\vspace{0.3cm}                                
$8$  & $5$  & $11$ & $10$ \\                                
$12$ & $9$  & $7$  & $6$ 
\end{tabular}
$$

\vspace{0.3cm}
21. Plusieurs autres figures de ce genre peuvent être formées, dont le 
carré suivant : 
$$
\begin{tabular}{cccc}
\vspace{0.4cm}                                
$a\alpha$ & $d\beta$  & $a\delta$ & $d\gamma$ \\
\vspace{0.4cm}                                
$b\delta$ & $c\gamma$ & $b\alpha$ & $c\beta$  \\
\vspace{0.4cm}                                
$d\alpha$ & $a\beta$  & $d\delta$ & $a\gamma$ \\                         
$c\delta$ & $b\gamma$ & $c\alpha$ & $b\beta$ 
\end{tabular}
$$
où il est clair que pour les lettres latines il faut prendre $$a+d=b+c\,,$$ et
pour les lettres grecques $$\alpha+\delta=\beta+\gamma\,,$$ d'où, si nous 
prenons les mêmes valeurs qu'avant, apparaît le carré magique suivant :
$$
\begin{tabular}{cccc}
\vspace{0.3cm}                                
$1$  & $14$ & $4$  & $15$ \\
\vspace{0.3cm}                                
$8$  & $11$ & $5$  & $10$ \\
\vspace{0.3cm}                                  
$13$ & $2$  & $16$ & $3$  \\                              
$12$ & $7$  & $9$  & $6$ 
\end{tabular}
$$

\vspace{0.3cm}
22. Dans toutes ces dispositions, tant les lettres latines que grecques font la
même somme dans toutes les rangées ; il peut également arriver que ceci 
se produise même sans cet usage, bien que la somme de toutes [les rangées]
conserve la somme due. Mais il serait inutile de passer plus en revue les 
irrégularités de ce genre, puisque pour de tels cas on ne peut donner 
aucune règle sûre. C'est pourquoi parmi les catégories suivantes nous 
étudierons de préférence les cas où la signification des lettres tant 
latines que grecques n'est soumise à aucune restriction.
\end{quotation}

\subsection{La catégorie des carrés divisés en 25 cases}
\label{25}
\begin{quotation}

23. Ici apparaissent donc cinq lettres latines $a$, $b$, $c$, $d$, $e$ et cinq
grecques $\alpha$, $\beta$, $\gamma$, $\delta$, $\varepsilon$, dont les valeurs
seront respectivement 0, 5, 10, 15, 20, et 1, 2, 3, 4,~5 ; il faut donc inscrire
dans les cases du carré des couples [mixtes] de ces lettres en telle sorte que
l'on retrouve dans chaque rangée, tant horizontale que verticale, et même 
diagonale, chacune des lettres.

\vspace{0.3cm}
24. Inscrivons d'abord les lettres latines, dans l'ordre, dans la ligne
supérieure de ce carré, puis complétons la diagonale de gauche en telle 
sorte que dans aucune des rangées restantes la même lettre n'apparaisse deux
fois, ce qui peut être fait de plus d'une manière. Une fois que cette 
rangée a été constituée, remplissons l'autre diagonale de façon 
naturelle%
\footnote{La démarche est exactement la même que pour les carrés de 4
(voir \S.14).}, 
comme on peut le voir dans la figure ci-après :
$$
\begin{tabular}{ccccc}
\vspace{0.4cm}                                
$a\varepsilon$ & $b\delta$    & $c\gamma$    & $d\beta$    & $e\alpha$   \\
\vspace{0.4cm}                                
$e\beta$       & $c\alpha$    & $d\delta$    & $a\gamma$ & $b\varepsilon$ \\
\vspace{0.4cm}                                
$d\alpha$      & $e\gamma$    & $b\beta$  & $c\varepsilon$ & $a\delta$ \\  
\vspace{0.4cm}                               
$b\gamma$  & $d\varepsilon$   & $a\alpha$    & $e\delta$   & $c\beta$ \\
$c\delta$      & $a\beta$  & $e\varepsilon$  & $b\alpha$   & $d\gamma$
\end{tabular}
$$ 
Ensuite sous la case centrale il faut inscrire le $a$ et, au-dessus, le $d$,
ce qui fait que la colonne du milieu sera déjà complétée, et ensuite les
rangées restantes se remplissent naturellement.

\vspace{0.3cm}
25. Pour les lettres grecques, il n'est pas nécessaire de recourir à la
rangée diagonale, mais si nous examinons attentivement la colonne médiane, 
nous découvrons dans les cases correspondantes de part et d'autre [de 
celle-ci] deux lettres différentes ; c'est pourquoi dans cette rangée nous 
associons aux lettres latines individuelles les lettres grecques de même nom 
et aux endroits correspondants nous échangeons ces lettres grecques%
\footnote{Donc symétriquement par rapport à la colonne centrale.},
comme nous l'avons fait dans la figure [ci-dessus].

\vspace{0.3cm}
26. Dans cette figure, clairement, aucune limitation n'est prescrite, et de plus 
tant pour les lettres latines que grecques il est permis de prendre n'importe
quelles valeurs correspondantes ; par quoi, avec cinq lettres, pouvant former 
120 combinaisons, peuvent apparaître en tout 14\,400 configurations 
[différentes]%
\footnote{$5!=120$ combinaisons, et donc $120\cdot 120$ dispositions possibles (voir note p.\pageref{note}).}.

\vspace{0.3cm}
27. Si de plus voulons échanger entre elles ces lignes ou colonnes, nous
obtiendrons plusieurs autres formes, mais la plupart du temps celles-ci
nécessiteront certaines limitations à cause des diagonales ; par 
exemple si nous déplaçons la première colonne à la fin, apparaît
la forme suivante :
$$
\begin{tabular}{ccccc}
\vspace{0.4cm}                                
$b\delta$ & $c\gamma$ & $d\beta$ & $e\alpha$ & $a\varepsilon$ \\
\vspace{0.4cm}                                
$c\alpha$ & $d\delta$ & $a\gamma$ & $b\varepsilon$ & $e\beta$ \\
\vspace{0.4cm}                                
$e\gamma$ & $b\beta$ & $c\varepsilon$ & $a\delta$ & $d\alpha$ \\  
\vspace{0.4cm}                               
$d\varepsilon$ & $a\alpha$ & $e\delta$ & $c\beta$ & $b\gamma$ \\
\vspace{0.4cm}
$a\beta$ & $e\varepsilon$ & $b\alpha$ & $d\gamma$ & $c\delta$ 
\end{tabular}
$$ 
où chacune des lettres apparaît dans chaque ligne et colonne ; 
mais pour qu'en même temps la condition soit satisfaite dans les diagonales,
[il faut que] tant cette somme
$$
3c+b+d+3\delta+\beta+\varepsilon
$$
que celle-ci
$$
3a+b+c+3\varepsilon+\alpha+\beta
$$
fassent la somme prescrite de toutes les lettres latines et grecques, à savoir 
$$
a+b+c+d+e+\alpha+\beta+\gamma+\delta+\varepsilon\,,
$$
[et] puisque de ce qui précède on peut en tirer ces deux équations :
$$
2c+2\delta=a+e+\alpha+\gamma
$$
et  
$$
2a+2\varepsilon=d+e+\gamma+\delta\,,
$$
ces conditions pourront être satisfaites de plusieurs façons ; bien plus,
tant les lettres latines que les grecques pourront ainsi être déterminées
séparément%
\footnote{Dans le texte : \textit{seorsim}, au lieu de \textit{seorsum (adv.)}, mais ceci ne semble pas être une exception.}%
\footnote{De la même manière que l'on peut séparer partie réelle et
imaginaire d'un nombre complexe.},
en sorte d'avoir
$$
\begin{tabular}{rlrlrllrl}
1) & $\!\!2c=a+e$, & 2) & $\!\!2a=d+e$, & 3) & $\!\!2\delta=\alpha+\gamma$ & et\ \ \ 
4) & $\!\!2\varepsilon=\gamma+\delta\,.$
\end{tabular}
$$
Il est en effet évident que les deux premières [équations] ci-dessus 
seront satisfaites si les lettres $d$, $b$, $a$, $c$, $e$ forment une 
progression arithmétique, ce qui se fait en prenant 
$$
d=0,\ b=5,\ a=10,\ c=15\ \mbox{et } e=20\,;
$$
les deux conditions restantes seront remplies si les lettres grecques,
disposées dans l'ordre $\alpha$, $\beta$, $\delta$, $\varepsilon$, $\gamma$
apparaissent [également] en progression arithmétique, ce qui sera le cas en
prenant
$$
\alpha=1,\ \beta=2,\ \delta=3,\ \varepsilon=4\ \mbox{et } \gamma=5\,,
$$
donc il en résulte le carré :
$$
\begin{tabular}{ccccc}
\vspace{0.3cm}                                
$8$  & $20$ & $2$  & $21$ & $14$ \\
\vspace{0.3cm}                                
$16$ & $3$  & $15$ & $9$  & $22$ \\
\vspace{0.3cm} 
$25$ & $7$  & $19$ & $13$ & $1$  \\  
\vspace{0.3cm}            
$4$  & $11$ & $23$ & $17$ & $10$ \\
$12$ & $24$ &  $6$ &  $5$ & $18$
\end{tabular}
$$
autrement dit, on a [un carré qui a] partout la même somme = 65.

\vspace{0.3cm}
28. Mais une telle attribution des lettres n'est pas une mince affaire et 
demande de l'attention, surtout pour les carrés d'ordres plus élevés, où
une plus grande quantité de lettres est laissée à notre bon vouloir%
\footnote{Donc qu'il nous revient de les placer selon notre choix.},
en sorte que le nombre de telles figures devient toujours plus grand ; 
mais si nous voulons laisser de côté cette condition, comme elle n'est 
soumise à aucune restriction sur les valeurs des lettres, le travail peut 
être rendu assez agréable ; si en effet la valeur du milieu [de la 
progression arithmétique], c'est-à-dire 10, est attribuée à la lettre 
$c$, laissant les autres valeurs à notre choix arbitraire, nous pourrons 
remplir l'une des deux diagonales avec cette même lettre $c$, et mettre 
les lettres restantes dans l'ordre alphabétique, de la manière que l'on voit
clairement ressortir de cette figure :
$$
\begin{tabular}{ccccc}
\vspace{0.1cm}                                
$c$ & $d$ & $e$ & $a$ & $b$ \\
\vspace{0.1cm}                                
$b$ & $c$ & $d$ & $e$ & $a$ \\
\vspace{0.1cm} 
$a$ & $b$ & $c$ & $d$ & $e$ \\
\vspace{0.1cm}            
$e$ & $a$ & $b$ & $c$ & $d$ \\
$d$ & $e$ & $a$ & $b$ & $c$ 
\end{tabular}
$$
Maintenant associons aux lettres latines de la ligne médiane les lettres
grecques de même nom, et ensuite échangeons deux à deux autour de cette
ligne les lettres grecques de même nom ; de cette opération apparaît
la forme suivante :
$$
\begin{tabular}{ccccc}
\vspace{0.4cm} 
$c\delta$ & $d\varepsilon$ & $e\alpha$ & $a\beta$ & $b\gamma$ \\
\vspace{0.4cm} 
$b\varepsilon$ & $c\alpha$ & $d\beta$ & $e\gamma$ & $a\delta$ \\
\vspace{0.4cm} 
$a\alpha$ & $b\beta$ & $c\gamma$ & $d\delta$ & $e\varepsilon$ \\
\vspace{0.4cm} 
$e\beta$ & $a\gamma$ & $b\delta$ & $c\varepsilon$ & $d\alpha$ \\
$d\gamma$ & $e\delta$ & $a\varepsilon$ & $b\alpha$ & $c\beta$ 
\end{tabular}
$$
d'où il est manifeste qu'il faut prendre pour $\gamma$ la valeur médiane, 
qui est 3 ; donc si nous posons dans l'ordre
$$
a=0,\ b=5,\ c=10,\ d=15,\ e=20
$$
et
$$
\alpha=1,\ \beta=2,\ \gamma=3,\ \delta=4,\ \varepsilon=5\,.
$$
il en résulte le carré magique suivant :
$$
\begin{tabular}{ccccc}
\vspace{0.3cm}                                
$14$ & $20$ & $21$ & $2$ & $8$ \\
\vspace{0.3cm}                                
$10$ & $11$ & $17$ & $23$ & $4$ \\
\vspace{0.3cm}                                
$1$ & $7$ & $13$ & $19$ & $25$ \\
\vspace{0.3cm}                                
$22$ & $3$ & $9$ & $15$ & $16$ \\                       
$18$ & $24$ & $5$ & $6$ & $12$ 
\end{tabular}
$$

\vspace{0.3cm}
29. Au moyen de la règle usuelle concernant la formation des carrés
impairs%
\footnote{La règle de construction est la suivante (voir la figure) : on place le 1
dans la case directement sous la case centrale ; on place ensuite le 2 dans
la case en diagonale en bas à droite du 1, et ainsi de suite dans l'ordre 
croissant (en suivant les diagonales brisées), jusqu'à avoir placé $n$ nombres ($n$ étant l'ordre du carré considéré ; on se retrouve ainsi diagonalement en haut à gauche du 1) ; on descend ensuite de deux cases vers le bas, et ainsi de suite jusqu'à ce que le carré soit rempli. Successivement appelée \textit{méthode de Bachet} (\textit{Problemes plaisans et delectables qui se font par les nombres}, 1612), puis \textit{de Cardan} (1550), puis \textit{de Moschopoulos} (env. 1300, Byzance), cette méthode est justifiée au début du XI$^{\grave{e}me}$ siècle chez Ibn al-Haytham (env. 965-1040). Son texte est perdu, mais la méthode nous est rapportée par un auteur anonyme du XII$^{\grave{e}me}$ siècle.},
qu'il est coutume de rappeler partout, on forme cette figure
$$
\begin{tabular}{ccccc}
\vspace{0.3cm}                                
$11$ & $24$ & $7$ & $20$ & $3$ \\
\vspace{0.3cm}                                
$4$ & $12$ & $25$ & $8$ & $16$ \\
\vspace{0.3cm}                                
$17$ & $5$ & $13$ & $21$ & $9$ \\
\vspace{0.3cm}                                
$10$ & $18$ & $1$ & $14$ & $22$ \\            
$23$ & $6$ & $19$ & $2$ & $15$ 
\end{tabular}
$$
dont on peut se demander si elle fait partie de notre catégorie ; observons 
tout d'abord que pour la diagonale de gauche, à cause de $c=10$, il faut poser 
$$
\delta=1,\ \alpha=2,\ \gamma=3,\ \varepsilon=4\ \mbox{et } \beta=5,
$$
et donc
$$
b=0,\ d=20,\ a=15,\ e=5\,.
$$
De ces valeurs il naît précisément ce carré.

\vspace{0.3cm}
30. On peut découvrir plusieurs autres formes fort régulières de cette
sorte%
\footnote{Allusion au mode de placement précédent.}, 
tant pour cette espèce que pour les suivantes ; ce faisant le nombre de 
carrés magiques pourra sans peine être augmenté d'une manière
considérable. Nous ne serions nullement sûrs d'avoir un jour épuisé tous
les cas possibles, même si leur nombre n'est certainement pas infini. Et il 
serait sans doute très souhaitable que l'on trouve des règles plus 
générales et appropriées à l'usage pratique, pour ne pas avoir besoin
d'opérations utilisant le tâtonnement. En effet on parviendrait ainsi à
une avancée fort belle dans le domaine de la théorie des combinaisons.
\end{quotation}

\subsection{La catégorie des carrés divisés en 36 cases}
\label{36}
\begin{quotation}

31. Puisqu'ici le nombre de formes différentes est très grand et que de
nombreuses déterminations sont laissées à notre bon jugement, apportons
ici seulement une règle particulière au moyen de laquelle les lettres, tant
latines que grecques, peuvent être facilement placées dans l'ordre dû, à
savoir qu'il faut attribuer aux six lettres latines des valeurs telles que l'on
ait
$$
a+f=b+e=c+d
$$
et de la même manière pour les lettres grecques
$$
\alpha+\zeta=\beta+\varepsilon=\gamma+\delta\,;
$$
alors en effet, par analogie avec le \S.20, inscrivons dans chaque ligne 
les deux lettres latines conjuguées ; quant aux colonnes, disposons-y de la 
même manière les paires de lettres grecques. De cette façon, on 
obtiendra la figure suivante%
\footnote{Ceci toutefois n'aboutit pas à un carré magique : en effet, dans
l'une des diagonales apparaît deux fois $b\beta$, dans l'autre deux fois
$e\varepsilon$.\ \ \ L.G.D. (Cette remarque laisse penser qu'il manque quelque 
chose après la figure. De plus, la méthode de la rangée médiane ne peut
pas s'appliquer puisque l'on a deux fois $e$ de part et d'autre de la
diagonale de gauche, alors qu'il est exigé d'avoir des lettres deux à deux
différentes (voir \S.14 et \S.24).)} :
$$
\begin{tabular}{cccccc}
\vspace{0.4cm} 
$a\alpha$ & $a\zeta$ & $a\beta$ & $f\varepsilon$ & $f\gamma$ & $f\delta$ \\
\vspace{0.4cm} 
$f\alpha$ & $f\zeta$ & $f\beta$ & $a\varepsilon$ & $a\gamma$ & $a\delta$ \\
\vspace{0.4cm} 
$b\alpha$ & $b\zeta$ & $b\beta$ & $e\varepsilon$ & $e\gamma$ & $e\delta$ \\
\vspace{0.4cm} 
$e\zeta$ & $e\alpha$ & $e\varepsilon$ & $b\beta$ & $b\delta$ & $b\gamma$ \\
\vspace{0.4cm} 
$c\zeta$ & $c\alpha$ & $c\varepsilon$ & $d\beta$ & $d\delta$ & $d\gamma$ \\ 
$d\zeta$ & $d\alpha$ & $d\varepsilon$ & $c\beta$ & $c\delta$ & $c\gamma$ 
\end{tabular}
$$

\vspace{0.3cm}
32. De là l'on comprend suffisamment clairement qu'une telle disposition
pourra s'appliquer avec succès à tous les carrés pairs%
\footnote{Mais voyez la note précédente et aussi les Commentaires 530 de ce
volume. Dans les \textit{Adversariis mathematicis}, dont les fragments sont 
édités à la page 535, on trouve ce carré magique de 36 cases
$$
\begin{tabular}{cccccc}
\vspace{0.3cm} 
$3$ & $36$ & $30$ & $4$ & $11$ & $27$ \\
\vspace{0.3cm} 
$22$ & $13$ & $35$ & $12$ & $14$ & $15$ \\
\vspace{0.3cm} 
$16$ & $18$ & $8$ & $31$ & $17$ & $21$ \\
\vspace{0.3cm} 
$28$ & $20$ & $6$ & $29$ & $19$ & $9$ \\
\vspace{0.3cm} 
$32$ & $23$ & $25$ & $2$ & $24$ & $5$ \\
$10$ & $1$ & $7$ & $33$ & $26$ & $34$ 
\end{tabular}
$$
qui, selon {\sc Euler} (voir p. 536) doit être appelé [carré] parfait,
parce que la même somme 111 se retrouve partout, tant dans les lignes et les
colonnes que dans les deux diagonales.\ \ \ L.G.D. 
(N.B. : ce carré ne s'y trouve pas d'ailleurs.)},
de même que la méthode décrite précédemment pour les carrés impairs
dans laquelle les lettres possédant les valeurs moyennes sont répétées
continûment, et ensuite les lettres restantes sont disposées à la suite
dans l'ordre naturel%
\footnote{Voir \S.28.},
en sorte que, quel que soit le nombre de cases proposé dans le carré, nous
ayons toujours la possibilité de construire beaucoup de carrés magiques,
même si les règles rapportées ici sont très particulières.
\end{quotation}

\pagebreak
\section{Recherche sur une nouvelle espèce de carrés magiques}
\begin{quotation}

1. Comme à son habitude, {\sc Euler} commence son mémoire par la description 
du contexte dans lequel il a été amené à étudier le problème des 36 
officiers. Il s'agit d'<<une question fort curieuse, qui a exercé pendant 
quelque temps la sagacité de bien du monde>> et qui l'a <<engagé à faire 
les recherches suivantes, qui semblent ouvrir une nouvelle carrière dans 
l'Analyse et en particulier dans la doctrine des combinaisons>>. 

Le problème est le suivant : il s'agit de 36 officiers, de six différents 
grades et tirés de six régiments différents, qu'il s'agit de ranger dans 
un carré de manière que sur chaque ligne et colonne, il se trouve six 
officiers tant de différents caractères que de régimes différents. 

{\sc Euler} précise immédiatement qu'un tel arrangement est impossible, 
malgré tous les efforts employés pour résoudre le problème, bien qu'il 
ne puisse pas en donner de démonstration rigoureuse.

\vspace{0.3cm}

2. Il commence par décrire la méthode qu'il a utilisée pour essayer 
de résoudre le problème : il marque les six régiments par les lettres 
latines $a$, $b$, $c$, $d$, $e$, $f$, et les six différents grades par les 
lettres grecques $\alpha$, $\beta$, $\gamma$, $\delta$, $\varepsilon$, $\zeta$ ;
il est assez clair que le caractère unique de chaque officier est 
entièrement déterminé par une combinaison de deux lettres, l'une latine et
l'autre grecque. 

Le problème peut donc être ainsi exprimé de la façon suivante : il 
s'agit d'inscrire ces 36 termes dans les 36 cases d'un carré, en sorte que 
l'on rencontre sur chaque ligne et colonne tant les six lettres latines que les 
grecques.

\vspace{0.3cm}

3. Il y a donc trois conditions à remplir :

\begin{description}
    \item[(a)] on doit retrouver sur chaque ligne les six lettres latines 
    et grecques ;
    \item[(b)] de même sur les colonnes ;
    \item[(c)] on doit retrouver effectivement tous les 36 termes dans le 
    carré (ou, ce qui revient au même, qu'aucun terme ne se retrouve 
    deux fois).
\end{description}

{\sc Euler} précise ici que les deux premières conditions ne suffisent pas, 
car sinon <<il ne seroit pas difficile de trouver plusieures solutions>>. Pour
confirmer ses dires, il donne l'exemple suivant :
$$
\begin{tabular}{cccccc}
\vspace{0.4cm}
$a\alpha$ & $b\zeta$ & $c\delta$ & $d\varepsilon$ & $e\gamma$ & $f\beta$ \\
\vspace{0.4cm}
$b\beta$ & $c\alpha$ & $f\varepsilon$ & $e\delta$ & $a\zeta$ & $d\gamma$ \\
\vspace{0.4cm}
$c\gamma$ & $d\varepsilon$ & $a\beta$ & $b\zeta$ & $f\delta$ & $e\alpha$ \\
\vspace{0.4cm}
$d\delta$ & $f\gamma$ & $e\zeta$ & $c\beta$ & $b\alpha$ & $a\varepsilon$ \\
\vspace{0.4cm}
$e\varepsilon$ & $a\delta$ & $b\gamma$ & $f\alpha$ & $d\beta$ & $c\zeta$ \\
$f\zeta$ & $e\beta$ & $d\alpha$ & $a\gamma$ & $c\varepsilon$ & $b\delta$ 
\end{tabular}
$$
où il constate que les termes $b\zeta$ et $d\varepsilon$ s'y rencontrent deux 
fois et que les termes $b\varepsilon$ et $d\zeta$ n'y sont pas.

\vspace{0.3cm}

4. {\sc Euler} s'avoue vaincu de n'avoir pu trouver de solution à la 
construction d'un tel carré, et, <<pour donner plus d'étendue à ses 
recherches>>, il passe à la généralisation du problème en 
considérant cette fois $n$ grades, symbolisés par les lettres latines $a$, 
$b$, $c$, $d$, etc., et autant de régiments, symbolisés par les lettres 
grecques $\alpha$, $\beta$, $\gamma$, $\delta$, etc. à combiner de $n^2$ 
manières différentes dans un carré de $n^2$ cases, avec les mêmes trois 
conditions qu'auparavant%
\footnote{Généralisées à $n$ termes, bien évidemment.}.

\vspace{0.3cm}

5. {\sc Euler} remarque ensuite que, si l'on donne aux lettres latines les 
valeurs 0, $n$, $2n$, ..., $(n-1)n$ et aux lettres grecques les valeurs 1, 2, 3,
..., $n$, on peut ainsi produire tous les nombres dans l'ordre naturel%
\footnote{Voir le \S.3 de la section 2.}, 
et comme chaque ligne (et colonne) du carré doit contenir toutes les 
différentes lettres (et donc la somme sur chaque ligne, respectivement chaque 
colonne, sera la même), alors un tel arrangement doit satisfaire à la 
condition des carrés magiques ordinaires. A la différence près que ces 
derniers n'ont pas besoin d'avoir toutes les valeurs dans chaque rangée. 

\vspace{0.3cm}

6. {\sc Euler} désire rendre <<plus commode les opérations à faire dans la 
suite>>. Il décide d'adopter la notation suivante : il met en lieu et place 
des lettres latines et grecques les nombres naturels 1, 2, 3, etc., et pour les
distinguer entre eux, il les appelle les uns \textit{nombres latins} et les autres 
\textit{nombres grecs} ; de plus, pour ne pas les confondre, il joint les nombres 
grecs aux latins en forme d'exposants (ainsi par exemple, $1^3$ signifiera 
$a\gamma$, suivant la notation que nous avons vue à la section
précédente).

Il donne comme exemple le carré de 7 suivant : 
$$
\begin{tabular}{lllllll}
\vspace{0.3cm}
$1^1$ & $2^6$ & $3^4$ & $4^3$ & $5^7$ & $6^5$ & $7^2$ \\
\vspace{0.3cm}
$2^2$ & $3^7$ & $1^5$ & $5^4$ & $4^1$ & $7^6$ & $6^3$ \\
\vspace{0.3cm}
$3^3$ & $6^1$ & $5^6$ & $7^5$ & $1^2$ & $4^7$ & $2^4$ \\
\vspace{0.3cm}
$4^4$ & $5^2$ & $6^7$ & $1^6$ & $7^3$ & $2^1$ & $3^5$ \\
\vspace{0.3cm}
$5^5$ & $1^3$ & $7^1$ & $2^7$ & $6^4$ & $3^2$ & $4^6$ \\
\vspace{0.3cm}
$6^6$ & $7^4$ & $4^2$ & $3^1$ & $2^5$ & $5^3$ & $1^7$ \\
$7^7$ & $4^5$ & $2^3$ & $6^2$ & $3^6$ & $1^4$ & $5^1$ 
\end{tabular}
$$
où il a rangé les nombres latins dans l'ordre naturel dans la première ligne 
et la première colonne, et égalé les nombres grecs en exposant aux nombres 
latins dans la première colonne, <<comme partout dans la suite, puisque la 
signification de ces nombres est absolument arbitraire>>.

\vspace{0.3cm}

7. Le lecteur se convaincra assez vite que le carré ci-dessus satisfait à la 
fois aux trois conditions requises. 

{\sc Euler} explique ensuite la construction de ce carré : il part tout 
d'abord de ce qu'il appelle le \textit{carré latin} fondamental, qui possède 
tous les sept nombres dans chacune des lignes et colonnes, et qui a donc la 
forme suivante :
$$
\begin{tabular}{ccccccc}
\vspace{0.1cm}
$1$ & $2$ & $3$ & $4$ & $5$ & $6$ & $7$ \\
\vspace{0.1cm}
$2$ & $3$ & $1$ & $5$ & $4$ & $7$ & $6$ \\
\vspace{0.1cm}
$3$ & $6$ & $5$ & $7$ & $1$ & $4$ & $2$ \\
\vspace{0.1cm}
$4$ & $5$ & $6$ & $1$ & $7$ & $2$ & $3$ \\
\vspace{0.1cm}
$5$ & $1$ & $7$ & $2$ & $6$ & $3$ & $4$ \\
\vspace{0.1cm}
$6$ & $7$ & $4$ & $3$ & $2$ & $5$ & $1$ \\
$7$ & $4$ & $2$ & $6$ & $3$ & $1$ & $5$ 
\end{tabular}
$$

\vspace{0.3cm}

8. Il faut maintenant trouver une méthode sûre pour joindre les nombres grecs 
à chaque nombre latin. Chaque nombre grec doit donc se retrouver à des 
niveaux différents dans les lignes, et de même pour les colonnes. 

\vspace{0.3cm}

9. {\sc Euler} procède de la façon suivante : il appelle \textit{formules 
directrices} les formules qui lui permettent de régler l'inscription des 
exposants. Revenant à l'exemple du \S.6, il décrit les formules directrices 
pour chaque exposant :
\vspace{0.3cm}
\begin{center}
\begin{tabular}{c|ccccccc}
\mbox{exposant} & \multicolumn{7}{c}{formule directrice} \\
\hline
1 & 1 & 6 & 7 & 3 & 4 & 2 & 5 \\
2 & 2 & 5 & 4 & 6 & 1 & 3 & 7 \\ 
3 & 3 & 1 & 2 & 4 & 7 & 5 & 6 \\
4 & 4 & 7 & 3 & 5 & 6 & 1 & 2 \\
5 & 5 & 4 & 1 & 7 & 2 & 6 & 3 \\
6 & 6 & 2 & 5 & 1 & 3 & 7 & 4 \\
7 & 7 & 3 & 6 & 2 & 5 & 4 & 1
\end{tabular}
\end{center}
\vspace{0.3cm}
Donc pour chaque exposant il repère successivement leur base dans chaque 
colonne, et la suite de ces bases forme sa formule directrice.

{\sc Euler} mentionne ensuite que pour avoir un carré complet il faut avoir 
une formule directrice pour chaque exposant, et que ces formules <<s'accordent 
tellement entre elles qu'en les écrivant l'une sous l'autre on rencontre dans 
chaque rangée verticale tous les différents nombres>>. Ceci est clair : 
sinon l'un des nombres latins aurait alors deux fois le même exposant.

\vspace{0.3cm}

10. Donc une fois établi un carré latin d'ordre quelconque, le travail 
consiste à chercher les formules directrices pour chaque exposant. Le lecteur 
se convaincra aisément que si l'on n'arrive pas à trouver une telle formule,
ne serait-ce que pour un seul exposant, alors la construction du carré est
impossible, de même si ces formules directrices ne s'accordent pas entre
elles.

{\sc Euler} met toutefois en garde qu'il faut bien s'assurer d'avoir examiné
tous les cas possibles avant d'arriver à une telle conclusion.

\vspace{0.3cm}

11. La première méthode utilisée pour trouver les formules directrices
pour chaque exposant est le tâtonnement ! {\sc Euler} décrit sa recherche
pour l'exposant~4 de son exemple, en choisissant à volonté les quatre
premiers termes, à savoir
$$
\begin{array}{cccc}
4 & 7 & 3 & 5
\end{array}
$$
qui sont tirés des quatre premières colonnes et des lignes 1, 2, 4 et 6.

Les trois nombres restants étant 1, 2, et 6, il faut donc les trouver dans les
trois dernières colonnes et parmi les lignes 3, 5 et 7, c'est-à-dire dans
le \textit{morceau}%
\footnote{{\sc Euler} utilise ce terme pour désigner toute partie du carré.}
suivant :
$$
\begin{array}{ccc}
1 & 4 & 2 \\
6 & 3 & 4 \\
3 & 1 & 5 
\end{array}
$$
duquel on tire facilement nos trois termes restants répondant aux conditions
prescrites.

\vspace{0.3cm}

12. {\sc Euler} recherche ensuite des formes plus particulières de 
carrés, auxquels des méthodes de recherche plus spécifiques peuvent être
appliquées. C'est dans ce but qu'il introduit ce qu'il appelle les {\it
carrés latins%
\footnote{Ce terme a d'ailleurs acquis droit de cité en mathématique.}
 à marche simple} (respectivement double, triple et quadruple).

\vspace{0.3cm}

13. Le carré à simple marche est celui où tous les nombres 1, 2, 3, ...,
$n$ <<marchent dans leur ordre naturel>>. Il aura la forme suivante :
$$
\begin{array}{cccccc}
1 & 2 & 3 & 4 & ... & n \\
2 & 3 & 4 & ... & n & 1 \\
3 & 4 & ... & n & 1 & 2 \\
4 & ... & n & 1 & 2 & 3 \\
\vdots & \vdots & \vdots & \vdots & \vdots & \vdots 
\end{array}
$$

\vspace{0.3cm}

14. Le \textit{carré à double marche} contient les nombres dans leur ordre 
naturel dans la première ligne et colonne, et deux à deux dans les 
suivantes. Il aura donc la forme générale suivante :
$$
\begin{array}{ccccccccc}
1 & 2 & 3 & 4 & 5  & 6  & 7  & 8  & ... \\
2 & 1 & 4 & 3 & 6  & 5  & 8  & 7  & ... \\
3 & 4 & 5 & 6 & 7  & 8  & 9  & 10 & ... \\
4 & 3 & 6 & 5 & 8  & 7  & 10 & 9  & ... \\
5 & 6 & 7 & 8 & 9  & 10 & 11 & 12 & ... \\
6 & 5 & 8 & 7 & 10 & 9  & 12 & 11 & ... \\
\vdots & \vdots & \vdots & \vdots & \vdots & \vdots & \vdots & \vdots & \ddots
\end{array}
$$

On remarque entre autres que cette espèce de carré ne peut exister que pour
des carrés pairs (la raison est évidente).

\vspace{0.3cm}

15. Les \textit{carrés à triple marche} sont construits de la même manière,
et n'ont de raison d'être que pour des carrés d'ordre divisible par 3. La 
forme générale sera donc la suivante :
$$
\begin{array}{cccccccccc}
1 & 2 & 3 & 4  & 5  & 6  & 7  & 8  & 9  & ... \\
2 & 3 & 1 & 5  & 6  & 4  & 8  & 9  & 7  & ... \\
3 & 1 & 2 & 6  & 4  & 5  & 9  & 7  & 8  & ... \\
4 & 5 & 6 & 7  & 8  & 9  & 10 & 11 & 12 & ... \\
5 & 6 & 4 & 8  & 9  & 7  & 11 & 12 & 10 & ... \\
6 & 4 & 5 & 9  & 7  & 8  & 12 & 10 & 11 & ... \\
7 & 8 & 9 & 10 & 11 & 12 & 13 & 14 & 15 & ... \\
\vdots & \vdots & \vdots & \vdots & \vdots & \vdots & \vdots & \vdots & \vdots & \ddots
\end{array}
$$

Il faut encore remarquer que cette forme est la seule possible, car l'on ne peut
pas former d'autres combinaisons à partir du premier membre suivant :
$$
\begin{array}{ccc}
1 & 2 & 3 \\
2 & 3 & 1 \\
3 & 1 & 2 
\end{array}
$$ 

En effet, une fois la première ligne et colonnes placées, on est obligé
de prendre~3 dans la case du milieu, car sinon il apparaît sous le 3 de la
première ligne, ce qui est interdit. Le reste se remplit alors de lui-même.

\vspace{0.3cm}

16. Pour les \textit{carrés à quadruple marche} la situation est différente, 
car on parvient à former 4 sortes de premiers membres, à savoir 
$$
\begin{array}{cccc|cccc|cccc|cccc}
\multicolumn{4}{c}{\mbox{I}}   & \multicolumn{4}{c}{\mbox{II}} & 
\multicolumn{4}{c}{\mbox{III}} & \multicolumn{4}{c}{\mbox{IV}} \\ 
1 & 2 & 3 & 4 & 1 & 2 & 3 & 4 & 1 & 2 & 3 & 4 & 1 & 2 & 3 & 4 \\
2 & 1 & 4 & 3 & 2 & 1 & 4 & 3 & 2 & 3 & 4 & 1 & 2 & 4 & 1 & 3 \\
3 & 4 & 1 & 2 & 3 & 4 & 2 & 1 & 3 & 4 & 1 & 2 & 3 & 1 & 4 & 2 \\
4 & 3 & 2 & 1 & 4 & 3 & 1 & 2 & 4 & 1 & 2 & 3 & 4 & 3 & 2 & 1
\end{array}
$$
et donc par là quatre formes générales. Comme les précédents, ceci
n'est valable que pour des carrés d'ordre divisible par 4.

{\sc Euler} examine ensuite en détail chaque espèce rapportée ci-dessus,
pour en déduire les carrés complets. Nous nous limiterons à étudier la
démarche utilisée pour les carrés à simple marche.
\end{quotation}

\subsection{Des carrés latins à simple marche de la forme générale}
\begin{quotation}

$$
\begin{array}{cccccc}
1 & 2 & 3 & 4 & ... & n \\
2 & 3 & 4 & ... & n & 1 \\
3 & 4 & ... & n & 1 & 2 \\
4 & ... & n & 1 & 2 & 3 \\
\vdots & \vdots & \vdots & \vdots & \vdots & \vdots 
\end{array}
$$

\pagebreak

\begin{center}
CAS DE $n=2$
\end{center}

17. Le cas le plus simple est vite réglé : le carré latin est 
$$
\begin{array}{cc}
1 & 2 \\
2 & 1
\end{array}
$$
et il est évident qu'on ne peut en tirer aucune formule directrice%
\footnote{En effet, en suivant la démarche on obtiendrait le carré
$$
\begin{array}{cc}
\vspace{0.2cm}
1^1 & 2^2 \\
2^2 & 1^1
\end{array}
$$
et donc $1^1$ et $2^2$ se retrouvent deux fois, et $2^1$ et $1^2$ manquent
totalement.}.

\vspace{0.3cm}
\begin{center}
CAS DE $n=3$
\end{center}

18. Le carré latin a la forme suivante :
$$
\begin{array}{ccc}
1 & 2 & 3 \\
2 & 3 & 1 \\
3 & 1 & 2
\end{array}
$$
dont la diagonale fournit tout de suite une formule directrice pour 
l'exposant~1 ; les autres se trouvent alors immédiatement. Puisque la formule 
1, 3, 2 est la seule possible pour cet exposant, on obtient le carré complet%
\footnote{Une tel carré est appelé \textit{carré eulérien.}}
 unique
$$
\begin{array}{ccc}
\vspace{0.2cm}
1^1 & 2^3 & 3^2 \\
\vspace{0.2cm}
2^2 & 3^1 & 1^3 \\
3^3 & 1^2 & 2^1
\end{array}
$$

\vspace{0.3cm}
\begin{center}
CAS DE $n=4$
\end{center}

19. Le carré est de la forme :
$$
\begin{array}{cccc}
1 & 2 & 3 & 4 \\
2 & 3 & 4 & 1 \\
3 & 4 & 1 & 2 \\
4 & 1 & 2 & 3
\end{array}
$$
mais dans ce cas il est impossible de trouver quelque formule directrice, quel
que soit l'exposant. {\sc Euler} en conclut donc l'impossibilité de former un
carré complet à partir du carré latin de 4 à simple marche.

\vspace{0.3cm}

20. Comme le même inconvénient se retrouve dans tous les cas où $n$ est
pair, {\sc Euler} énonce et démontre le théorème suivant :

\pagebreak

\begin{thm}\end{thm}
Pour tous les cas où le nombre $n$ est pair, le quarré latin à simple
marche ne sauroit jamais fournir une solution de la question proposée.

\vspace{0.3cm}
\begin{center}
CAS DE $n=5$
\end{center}

21. Pour le carré de 5, qui est de la forme 
$$
\begin{array}{ccccc}
1 & 2 & 3 & 4 & 5 \\
2 & 3 & 4 & 5 & 1 \\
3 & 4 & 5 & 1 & 2 \\
4 & 5 & 1 & 2 & 3 \\
5 & 1 & 2 & 3 & 4 
\end{array}
$$
il existe seulement trois formules directrices pour l'exposant 1, à savoir
$$
\begin{array}{ccccc}
1 & 3 & 5 & 2 & 4 \\
1 & 4 & 2 & 5 & 3 \\
1 & 5 & 4 & 3 & 2 
\end{array}
$$

{\sc Euler} déclare :
<<En ajoutant l'unité à chacun des termes de ces directrices, on obtiendra
celles pour l'exposant 2, qui, en y ajoutant à nouveau l'unité, donneront
celles pour l'exposant 3, et ainsi des autres>>%
\footnote{Il énonce cette règle de formation d'autres directrices, une 
fois que l'une d'entre elles a été trouvée, lorsqu'il aborde le carré 
d'ordre 7 (voir \S.32).}.

Ayant donc trouvé les directrices pour les autres exposants à partir
des trois susmentionnées, on pourra construire trois carrés de directrices
<<propres à diriger l'inscription des exposans>>. Prenons à titre d'exemple
celui qui est obtenu par la première formule :
$$
\begin{array}{ccccc}
1 & 3 & 5 & 2 & 4 \\
2 & 4 & 1 & 3 & 5 \\
3 & 5 & 2 & 4 & 1 \\
4 & 1 & 3 & 5 & 2 \\
5 & 2 & 4 & 1 & 3 
\end{array}
$$

\vspace{0.3cm}

22. Le carré complet qui lui correspond est donc :
$$
\begin{array}{ccccc}
\vspace{0.2cm}
1^1 & 2^5 & 3^4 & 4^3 & 5^2 \\
\vspace{0.2cm}
2^2 & 3^1 & 4^5 & 5^4 & 1^3 \\
\vspace{0.2cm}
3^3 & 4^2 & 5^1 & 1^5 & 2^4 \\
\vspace{0.2cm}
4^4 & 5^3 & 1^2 & 2^1 & 3^5 \\
5^5 & 1^4 & 2^3 & 3^2 & 4^1 
\end{array}
$$

\vspace{0.3cm}

23. {\sc Euler} remarque encore les choses suivantes : 
\begin{itemize}
    \item les termes des formules directrices vont en progression
    arithmétique (de raison respectivement 2, 3 et 4)%
\footnote{Ici {\sc Euler} mentionne encore : <<...dans la quatrième de
    5 et ainsi des autres.>>, ce qui n'a pas de rapport avec le cas
    considéré, puisqu'il n'y a que trois directrices. Probablement
    a-t-il déjà en vue le cas général.} ;  
    \item les exposants de la première ligne pour les trois carrés 
    complets correspondent aux trois formules directrices%
\footnote{Cette observation sera elle aussi généralisée et l'on
    voit qu'elle facilitera grandement la construction du carré complet
    une fois qu'une directrice aura été trouvée.} ;
    \item après un échange de lignes, le premier carré complet 
    (trouvé au \S22) donne le carré <<très remarquable>> suivant :
    $$
    \begin{array}{ccccc}
    \vspace{0.2cm}
    1^1 & 2^5 & 3^4 & 4^3 & 5^2 \\
    \vspace{0.2cm}
    3^3 & 4^2 & 5^1 & 1^5 & 2^4 \\
    \vspace{0.2cm}
    5^5 & 1^4 & 2^3 & 3^2 & 4^1 \\
    \vspace{0.2cm}
    2^2 & 3^1 & 4^5 & 5^4 & 1^3 \\
    4^4 & 5^3 & 1^2 & 2^1 & 3^5 
    \end{array}
    $$
    qui contient les différentes lettres grecques et latines même dans
    les deux diagonales principales et les diagonales brisées%
\footnote{Un tel carré magique est par ailleurs appelé 
    \textit{pandiagonal}. Dès la fin du XIX$^{\grave{e}me}$ siècle, il
    reçoit même l'appellation de carré \textit{diabolique} ou encore 
    \textit{satanique} (dû au caractère <<surnaturel>> d'un tel carré, et
    à son utilisation dans les sciences occultes).} !
\end{itemize}

\vspace{0.3cm}
\begin{center}
CAS DE $n=7$
\end{center}

24. Les difficultés commencent. Le carré latin est le suivant :
$$
\begin{array}{ccccccc}
1 & 2 & 3 & 4 & 5 & 6 & 7 \\
2 & 3 & 4 & 5 & 6 & 7 & 1 \\
3 & 4 & 5 & 6 & 7 & 1 & 2 \\
4 & 5 & 6 & 7 & 1 & 2 & 3 \\
5 & 6 & 7 & 1 & 2 & 3 & 4 \\
6 & 7 & 1 & 2 & 3 & 4 & 5 \\
7 & 1 & 2 & 3 & 4 & 5 & 6
\end{array}
$$

{\sc Euler} trouve d'abord les cinq directrices croissantes (en progression
arithmé\-tique) pour l'exposant 1. Il en trouve ensuite 14 autres de
<<manière très embarassante expliquée au \S.8>>.

\vspace{0.3cm}

25. Il précise encore que <<le bel ordre qui règne dans les quarrés à
simple marche nous fournit des moyens très faciles pour trouver plusieures
telles formules, dès qu'on en a trouvé une seule>>. Il s'emploie donc
désormais à trouver des règles de formation d'autres directrices, une fois
que l'une d'entre elles a été trouvée.

\vspace{0.3cm}

26. A cette fin, il propose la démarche suivante : soit une
directrice quelconque se rapportant à l'exposant 1, et dont l'indice d'un
terme quelconque $x$ soit égal à~$t$%
\footnote{C'est-à-dire le numéro de sa colonne.}
; on a donc la situation suivante :
$$
\begin {array}{c|cccccccc}
t & 1 & 2 & 3 & 4 & 5 & 6 & ... & \\
\cline{1-8}
x & 1 & a & b & c & d & e & ... & \longleftarrow \mbox{formule directrice}
\end{array}
$$
et il faut dès lors remarquer :
\begin{itemize}
    \item qu'il faut trouver une transformation telle que 
    \begin{itemize} 
        \item pour $t=1$, on ait $x=1$ ;
        \item lorsqu'on attribue à $t$ toutes les valeurs de 1 
        jusqu'à $n$, $x$ prenne également toutes ces valeurs ;
        \end{itemize}
    \item puisque $t$ est le numéro de la colonne de laquelle $x$ est
    tiré, le numéro de la ligne de laquelle $x$ est tiré est
    $x-t+1$ par construction ; par conséquent, puisque les nombres
    $a$, $b$, $c$, etc. doivent être pris à chaque fois d'une ligne
    différente, alors $x-t+1$, et donc $x-t$, doit prendre toutes les
    valeurs de 1 jusqu'à $n$.
\end{itemize}

\vspace{0.3cm}

27. Ayant remarqué cela, {\sc Euler} énonce une première règle pour la
formation d'une nouvelle directrice
$$
\begin{array}{ccccccc}
1 & A & B & C & D & E & ...
\end{array}
$$
dont l'indice d'un terme quelconque $X$ est égal à $T$ :

\begin{regle}\end{regle}
Prenés $x$ pour l'indice et $t$ pour le terme qui lui répond.

\vspace{0.3cm}

En d'autres termes il faut poser
$$
\begin{array}{ccc}
T=x & \mbox{et} & X=t\,,
\end{array}
$$
et l'on vérifie aisément que toutes les conditions requises sont bien
remplies, à savoir 
\begin{itemize}
    \item qu'en donnant à $T$ toutes les valeurs possibles, non seulement 
    $X$ mais aussi $X-T$ prendront toutes les valeurs de 1 à $n$ ;
    \item si $T=1$, on aura bien $X=1$.
\end{itemize}

\vspace{0.3cm}

28. La deuxième règle énoncée par {\sc Euler} est la suivante :

\begin{regle}\end{regle}
On obtiendra toujours une autre nouvelle directrice en prenant
$$
\begin{array}{ccc}
T=t & \mbox{et} & X=1+t-x\,.
\end{array}
$$

\vspace{0.3cm}

29. En combinant ces deux règles, on peut trouver facilement des nouvelles
formules, qu'on pourra représenter de la manière suivante :
$$
\begin{array}{|c|c|c|c|c|c|}
\hline
      & \mbox{I} & \mbox{II} & \mbox{III} & \mbox{IV} & \mbox{V} \\
\hline
T=t   &     x    &     t     &   1+t-x    &     x     &  1+x-t   \\
X=x   &     t    &   1+t-x   &     t      &   1+x-t   &    x     \\
\hline
\hline
\mbox{VI} & \mbox{VII} & \mbox{VIII} & \mbox{IX} & \mbox{X} & \mbox{XI} \\
\hline
  1+t-x   &    2-x     &    1+x-t    &    2-t    &    2-x   &  2-t  \\
   2-x    &   1+t-x    &     2-t     &   1+x-t   &    2-t   &  2-x  \\
\hline
\end{array}
$$

Les cinq premières sont triviales par l'énoncé même des deux règles
principales. Pour ce qui est des six dernières, remarquons que si nous 
appliquons la deuxième règle à la règle V, nous obtenons la règle 
VII :
$$
\left\{
\begin{array}{rclclcc}
\tilde{T} & = & 1+X-T & = & 1+(1+t-x)-t & = & 2-x \\
\tilde{X} & = & X     & = & 1+t-x       &&
\end{array}
\right.
$$     
et ainsi de suite. Les six dernières règles mentionnées ne seront utiles 
qu'à partir
de l'ordre $n=9$, mais {\sc Euler} les mentionne déjà%
\footnote{Ce sera d'ailleurs sa
troisième règle principale de formation de directrices (voir \S.41).}.

Toutes ces règles vérifient les conditions requises.

\vspace{0.3cm}

30. A titre d'exemple, et pour nous familiariser un peu avec ces différentes
opérations, prenons par exemple la directrice 
$$
\begin{array}{ccccccc}
1 & 4 & 2 & 7 & 6 & 3 & 5
\end{array}
$$
et appliquons-lui la \textit{première règle} : le premier élément 
reste le 1 ;
pour $t=2$, nous avons $x=4$, donc dans la nouvelle directrice, nous aurons
$X=2$ pour $T=4$, c'est-à-dire que nous aurons le 2 à la quatrième place. 
Il en va de même pour les nombres restants et nous trouvons la directrice
$$
\begin{array}{ccccccc}
1 & 3 & 6 & 2 & 7 & 5 & 4\,.
\end{array}
$$
Pour l'application de la deuxième règle, le procédé est absolument
identique.

\vspace{0.3cm}

31. On trouve ainsi onze nouvelles directrices à partir d'une directrice 
quelconque, par simple utilisation du tableau.

\vspace{0.3cm}

32. De chacune des directrices trouvées on peut former un carré complet. Et 
c'est une des remarques précédentes qui va nous amener à l'une des 
formations générales de carré de directrices : en prenant une directrice 
quelconque 1, $a$, $b$, $c$, $d$, $e$, $f$, et en {\bf <<continuant ces nombres 
suivant leur ordre naturel, on aura les directrices pour les exposants 
suivants>>}. On obtient ainsi le carré de directrices suivant :
$$
\begin{array}{ccccccc}
\vspace{0.2cm}
1 &  a  &  b  &  c  &  d  &  e  &  f  \\ 
\vspace{0.2cm}
2 & a+1 & b+1 & c+1 & d+1 & e+1 & f+1 \\
\vspace{0.2cm}
3 & a+2 & b+2 & c+2 & d+2 & e+2 & f+2 \\
\vspace{0.2cm}
4 & a+3 & b+3 & c+3 & d+3 & e+3 & f+3 \\
\vspace{0.2cm}
5 & a+4 & b+4 & c+4 & d+4 & e+4 & f+4 \\
\vspace{0.2cm}
6 & a+5 & b+5 & c+5 & d+5 & e+5 & f+5 \\ 
7 & a+6 & b+6 & c+6 & d+6 & e+6 & f+6
\end{array}
$$
où il est assez clair que toutes les lignes et colonnes contiennent tous les
nombres de 1 à 7, quel que soit l'ordre des nombres $a$, $b$, $c$, $d$, $e$ et
$f$.

\vspace{0.3cm}

33. {\sc Euler} montre ensuite que si la directrice trouvée est 
1, $a$, $b$, $c$, $d$, $e$, $f$, la formule 

\begin{center}
1, $3-a$, $4-b$, $5-c$, $6-d$, $7-e$, $8-f$ 
\end{center}
est aussi une formule directrice pour la première ligne.

\vspace{0.3cm}

34. Comme nous avons déjà vu%
\footnote{Voir la remarque au \S.23.} que l'on pouvait prendre
comme exposants pour la première ligne la directrice quelconque trouvée, 
le carré complet qui apparaît alors sera donc :
$$
\begin{array}{ccccccc}
\vspace{0.2cm}
1^1 &   2^a   &   3^b   &   4^c   &   5^d   &   6^e   &   7^f   \\ 
\vspace{0.2cm}
2^2 & 3^{a+1} & 4^{b+1} & 5^{c+1} & 6^{d+1} & 7^{e+1} & 1^{f+1} \\
\vspace{0.2cm}
3^3 & 4^{a+2} & 5^{b+2} & 6^{c+2} & 7^{d+2} & 1^{e+2} & 2^{f+2} \\
\vspace{0.2cm}
4^4 & 5^{a+3} & 6^{b+3} & 7^{c+3} & 1^{d+3} & 2^{e+3} & 3^{f+3} \\
\vspace{0.2cm}
5^5 & 6^{a+4} & 7^{b+4} & 1^{c+4} & 2^{d+4} & 3^{e+4} & 4^{f+4} \\
\vspace{0.2cm}
6^6 & 7^{a+5} & 1^{b+5} & 2^{c+5} & 3^{d+5} & 4^{e+5} & 5^{f+5} \\ 
7^7 & 1^{a+6} & 2^{b+6} & 3^{c+6} & 4^{d+6} & 5^{e+6} & 6^{f+6}
\end{array}
$$

\vspace{0.3cm}

35. {\sc Euler} considère encore des carrés complets formés de plusieurs
formules directrices jointes, et en donne un exemple. 

\vspace{0.3cm}

36. Il s'étonne encore du fait que le cas $n=7$ fournisse autant de solutions,
alors que le cas $n=6$ n'en donne aucune et le cas $n=5$ trois seulement.

\vspace{0.3cm}
\begin{center}
CAS DE $n=9$
\end{center}

37-39. Puisque le nombre de directrices est énorme, {\sc Euler} se contente de
développer celles qui suivent une progression arithmétique, en excluant 
celles
dont la différence serait 3 ou 6 (car 3 et 6 ne sont pas premiers avec l'ordre
considéré). C'est là une remarque d'importance, car effectivement on se
convainc assez rapidement que si la différence des progressions possède un 
diviseur commun avec l'ordre~$n$, alors la formule directrice ne contiendra plus
tous les nombres de 1 jusqu'à $n$.

Excluant ces cas, les formules sont 
$$
\begin{tabular}{ccccccccc}
1 & 3 & 5 & 7 & 9 & 2 & 4 & 6 & 8 \\
1 & 6 & 2 & 7 & 3 & 8 & 4 & 9 & 5 \\
1 & 9 & 8 & 7 & 6 & 5 & 4 & 3 & 2
\end{tabular}
$$
On forme facilement les trois carrés complets suivant la même règle
que pour le cas précédent (\S.34).

\vspace{0.3cm}

40. {\sc Euler} considère encore les onze autres directrices obtenues en
appliquant les deux premières règles principales.

\vspace{0.3cm}

41. Il énonce ensuite sa troisième règle principale, dont il n'avait pas
fait usage jusqu'ici et que nous avons déjà vue au \S.29 : 

\begin{regle}\end{regle} Posant pour la directrice proposée l'indice $=t$ et
le terme qui lui répond $=x$, on pourra prendre pour la nouvelle formule
l'indice $T=2t-1$ et le terme même $X=2x-1$.

L'argument est le même que pour les deux premières règles : 
\begin{itemize}
    \item en posant $t=1$ et $x=1$ nous avons bien $T=1$ et $X=1$ ;
    \item donnant à $x$ toutes les valeurs de 1 à $n$, alors $2x$
    aussi (puisque l'ordre est impair) et donc $2x-1$ également ;
    \item si $x-1$ passe par toutes les valeurs, $X-T=2(x-t)$ prendra
    également toutes les valeurs possibles.
\end{itemize}

\vspace{0.3cm}

42. {\sc Euler} s'amuse ensuite à produire en grandes quantités des
directrices issues de cette nouvelle règle, et il cherche encore des 
directrices qui se reproduisent elle-mêmes par l'une des règles.

\vspace{0.3cm}

43-49. Ayant trouvé bon nombre de formules directrices (57 en tout), il en
expose encore huit autres trouvées par la <<méthode directe>>, et conclut
sans autre forme de procès que <<le nombre de toutes les directrices sera au 
moins quatre fois plus grand>>.
\end{quotation}

\subsection{Autres considérations}
\begin{quotation}

\vspace{0.3cm}

50-139. Dans la suite, {\sc Euler} considère encore les carrés dont les
diagonales principales et les diagonales brisées vérifient également les
conditions.

Il passe ensuite à l'étude des carrés latins à double, triple et
quadruple marche, dont il donne des règles de formation de directrices et 
plusieurs exemples.

\vspace{0.3cm}

140-152. Après avoir constaté qu'aucune des méthodes ne pouvait fournir de
solution au cas $n=6$ (qui, rappelons-le, est à la base de toute la
démarche) ainsi qu'aux cas où l'ordre est \textit{impairement pair}%
\footnote{C'est-à-dire de la forme $n=4k+2$, $k$=0, 1, 2, ...},
il suppose que si une telle solution existe, alors le carré doit être tout
à fait irrégulier, et qu'il faudrait examiner tous les cas possibles de
$n=6$.

Toute la cinquième partie de son travail est donc consacrée à la recherche
de carrés irréguliers. En particulier il donne une méthode <<par le moyen
de laquelle on peut transformer facilement, en plusieurs formes différentes,
tous les quarrés réguliers et examiner ensuite s'ils admettent des
directrices ou non>>.

La méthode tient dans le fait que l'on peut échanger entre eux deux nombres 
$a$ et $b$ si ceux-ci se trouvent dans les angles d'un rectangle, comme 
indiqué ci-dessous%
\footnote{La justification de cette façon de faire est évidente : les deux
lignes et colonnes concernées possèdent encore chacun des nombres après
l'échange.} :
$$
\begin{array}{ccc}
a      & ......... & b \\
\vdots &           & \vdots \\
b      & ......... & a
\end{array}
$$

A titre d'exemple, {\sc Euler} donne le carré à simple marche suivant :
$$
\begin{tabular}{cccccc}
\vspace{0.1cm}
1 & 2 & 3 & 4 & 5 & 6 \\
\vspace{0.1cm}
2 & \underline{3} & 4 & 5 & \underline{6} & 1 \\
\vspace{0.1cm}
3 & 4 & 5 & 6 & 1 & 2 \\
\vspace{0.1cm}
4 & 5 & 6 & 1 & 2 & 3 \\
\vspace{0.1cm}
5 & \underline{6} & 1 & 2 & \underline{3} & 4 \\
6 & 1 & 2 & 3 & 4 & 5
\end{tabular}
$$
dont il a montré au \S.20 qu'il n'avait aucune directrice. Or, en échangeant
les chiffres 3 et 6 de cet exemple, on obtient le carré suivant :
$$
\begin{tabular}{cccccc}
\vspace{0.1cm}
1 & 2 & 3 & 4 & 5 & 6 \\
\vspace{0.1cm}
2 & 6 & 4 & 5 & 3 & 1 \\
\vspace{0.1cm}
3 & 4 & 5 & 6 & 1 & 2 \\
\vspace{0.1cm}
4 & 5 & 6 & 1 & 2 & 3 \\
\vspace{0.1cm}
5 & 3 & 1 & 2 & 6 & 4 \\
6 & 1 & 2 & 3 & 4 & 5
\end{tabular}
$$
qui admet un grand nombre de directrices pour tous les six exposants
({\sc Euler} en trouve 32 !). 

Après avoir vérifié que pour chacune des directrices on ne pouvait
former aucun carré complet, et que ceci était le cas pour un très grand
nombre de carrés d'ordre~6, {\sc Euler} déclare : <<[...] je n'ai pas
hésité d'en conclure qu'on ne sauroit produire aucun quarré complet de 36
cases, et que la même impossibilité s'étende aux cas de $n=10$, $n=14$ et
en général à tous les nombres impairement pairs. [...] Or, ayant examiné
un nombre très considérable de tels quarrés, il me paroît impossible
que tous les cas mentionnés me fussent échappés.>>%
\footnote{Il est aujourd'hui prouvé qu'il existe des carrés latins pour n'importe quel ordre $n$, sauf $n=2$ et $n=6$.}

Cette certitude est encore accentuée avec l'étude d'une règle de
transformation générale : un carré peut être transformé en plusieurs
autres qui ont la même propriété par rapport aux directrices, si bien que
si le carré considéré n'en admet aucune, alors il en sera de même pour
tous les carrés transformés.

{\sc Euler} conclut en laissant le soin aux Géomètres%
\footnote{C'est-à-dire les mathématiciens.} de trouver des moyens
pour achever le dénombrement de tous les cas possibles, et se satisfait 
d'avoir apporté des observations <<assés importantes tant pour la doctrine 
des combinaisons que pour la théorie générale des quarrés magiques.>>
\end{quotation}

\pagebreak

\section{Conclusions}
\begin{quotation}
Le but principal de ce projet étant de revoir mes notions de latin, je
considère qu'il est partiellement atteint. Les nombreuses difficultés de la
traduction provenaient surtout de la grammaire, notamment les temps des verbes
et les formes conjugales de certains mots. 

Certains paragraphes sont très difficiles à traduire, alors que d'autres se
lisent pratiquement du français, ce qui a contribué a une avancée très
inégale et parfois très frustrante du travail. 

L'intérêt de traduire un texte mathématique est de pouvoir parfois
comprendre sa signification au moyen de l'exemple traité, car une formule
s'écrit évidemment de la même façon en latin qu'en français, et
ceci s'est avéré utile plus d'une fois.

En ce qui concerne la deuxième partie, je suis retombé en enfance le temps
d'un instant : une fois découvert qu'il s'agissait de l'analyse du
problème des 36 officiers, je me suis revu devant le tome 10 de ma collection
\textit{Alpha junior} en train de découvrir ce problème (vulgarisé) 
passionnant, et j'ai encore tête l'image d'un carré rempli 
d'officiers se regardant bizarrement...

\vspace{0.3cm}

Pour finir, j'ai remarqué lors de ce travail de semestre à quel point Dieu
avait guidé ma vie jusqu'ici : en effet, j'étais inscrit au gymnase en
section scientifique et la décision de finalement faire latin-langues est
intervenue juste avant de commencer les cours. De même, le fait d'avoir choisi
la section de mathématiques à l'EPFL, après plusieurs déboires en
microtechnique et en électricité, fut le résultat d'une collaboration
étroite avec mon Créateur.

Sans ce changement radical d'orientation, je n'aurais jamais eu la possibilité
de découvrir des textes mathématiques dans l'une des langues couramment 
utilisées à l'époque pour les rédiger, il est donc normal que je l'en
remercie. 

Je tiens encore à remercier M. Sesiano pour sa patience et ses
corrections nocturnes ($!$) de ma traduction et de ma grammaire <<germanisées>>,
ainsi que pour son enthousiasme lors de ses cours que j'ai eu le plaisir de 
suivre, car c'est à lui que je dois mon engouement pour l'histoire des
mathématiques.

\end{quotation}

\end{document}